\documentclass[a4paper,11pt]{amsart}
\usepackage[margin=32mm,bottom=40mm]{geometry}

\usepackage{lmodern}
\usepackage[stretch=10]{microtype}

\usepackage[T1]{fontenc}
\usepackage{stmaryrd, amsmath, amsthm, amssymb, amsfonts,
	hyperref, mathtools, mathrsfs, pythonhighlight}
\usepackage[shortlabels]{enumitem}
\usepackage{latexsym}
\usepackage{graphicx}
\usepackage{multicol,multirow}
\usepackage{amsmath,amssymb,amsfonts}
\usepackage{mathrsfs}
\usepackage{amsthm}
\usepackage{rotating}
\usepackage{appendix}
\usepackage[numbers]{natbib}
\usepackage[utf8x]{inputenc}

\numberwithin{equation}{section}

\newcommand{\bb}{\textbf}

\newcommand{\ol}{\overline}
\newcommand{\mc}{\mathcal}
\newcommand{\mf}{\mathfrak}

\newcommand{\II}{\mathbb{I}}
\newcommand{\ZZ}{\mathbb{Z}}
\newcommand{\CC}{\mathbb{C}}
\newcommand{\OO}{\mathbb{O}}
\newcommand{\RR}{\mathbb{R}}
\newcommand{\QQ}{\mathbb{Q}}

\newcommand{\FF}{\mathbb{F}}
\newcommand{\GG}{\mathbb{G}}
\DeclareMathOperator{\disc}{disc}
\DeclareMathOperator{\tr}{tr}
\DeclareMathOperator{\diag}{diag}

\DeclareMathOperator{\coker}{coker}
\DeclareMathOperator{\im}{im}
\DeclareMathOperator{\id}{id}
\DeclareMathOperator{\Span}{Span}
\DeclareMathOperator{\Gl}{Gl}
\DeclareMathOperator{\Sl}{Sl}
\DeclareMathOperator{\GU}{GU}
\DeclareMathOperator{\SU}{SU}
\DeclareMathOperator{\SO}{SO}
\DeclareMathOperator{\MT}{MT}
\DeclareMathOperator{\U}{U}
\DeclareMathOperator{\GCD}{GCD}
\DeclareMathOperator{\ord}{ord}
\DeclareMathOperator{\Spec}{Spec}
\DeclareMathOperator{\Gal}{Gal}
\DeclareMathOperator{\Hom}{Hom}
\DeclareMathOperator{\End}{End}

\DeclareMathOperator{\Lie}{Lie}
\DeclareMathOperator{\GSp}{GSp}
\DeclareMathOperator{\Sp}{Sp}
\DeclareMathOperator{\Hg}{Hg}

\DeclareMathOperator{\divv}{div}
\DeclareMathOperator{\Pic}{Pic}
\newcommand{\prs}[3]{{\left(\frac{#1}{#2} \right)_{#3}}}
\theoremstyle{plain}
\newtheorem{theorem}{Theorem}[section]
\newtheorem{lemma}[theorem]{Lemma}
\newtheorem{corollary}[theorem]{Corollary}
\newtheorem{proposition}[theorem]{Proposition}
\newtheorem{Setup}[theorem]{Setup}
\newtheorem{question}[theorem]{Question}

\theoremstyle{definition}

\newtheorem{remark}[theorem]{Remark}
\newtheorem{example}[theorem]{Example}

\begin{document}
\title{The ``exponential'' torsion of superelliptic Jacobians}
\author[J. Garnek]{J\k{e}drzej Garnek}
\address{\parbox{\linewidth}{Institute of Mathematics of Polish Academy of Sciences,\\ ul. \'{S}niadeckich 8, 00-656 Warszawa \\ \mbox{}}}
\email{jgarnek@amu.edu.pl}
\urladdr{http://jgarnek.faculty.wmi.amu.edu.pl/}
\date{}  
	
	\maketitle
	
	\begin{abstract}
		Let $J$ be the Jacobian of a superelliptic curve defined by the equation
		$y^{\ell} = f(x)$, where $f$ is a separable polynomial of degree non-divisible by $\ell$.
		Recall that $J$ is an abelian variety of dimension $\frac{1}{2}(\ell - 1)(\deg f-1)$.
		In this article we study the ``exponential'' (i.e. $\ell$-power) torsion of $J$. In particular, under some mild conditions on the polynomial~$f$, we determine
		the image of the associated $\ell$-adic representation up to the determinant. We show also that the image of the determinant is contained in an explicit $\ZZ_{\ell}$-lattice with a finite index.
		As an application, we prove new cases of the Hodge, Tate and Mumford--Tate conjectures for generic superelliptic Jacobians of the above type.
	\end{abstract}

%
\section{Introduction} \label{sec:intro}
\noindent 
A subject of a special study among abelian varieties are Jacobians of projective curves.
Often, when studying the Jacobians, one restricts further to some class of
curves with particularly nice properties, for example to superelliptic ones.
Recall that a \bb{superelliptic curve} over a field $K$ is a smooth projective curve $C$ with the affine part given by the equation of the form
\begin{equation} \label{eqn:superelliptic}
	C : y^m = f(x),
\end{equation}
where $f \in K[x]$ is separable of degree~$r$ and $m$ is a natural number (the \bb{exponent} of~$C$).
Similarly one says that $J$, the Jacobian of $C$, is superelliptic of exponent~$m$.
If $m = 2$, we call $C$ and $J$ hyperelliptic. 
There is a vast literature concerning  various properties of hyper- and superelliptic Jacobians, for instance their endomorphism rings (cf. \cite{Zarhin_Endomorphisms_ordinary_superelliptic},
\cite{Zarhin_endomorphisms_superelliptic}, \cite{Zarhin_Endomorphisms_reference_superelliptic}, \cite{Xue_Endomorphisms_superelliptic}, \cite{Zarhin_superelliptic_jacobians}),
their Hodge groups (cf. \cite{Xue_Zarhin_Hodge_groups},	\cite{Xue_Hodge_groups_superelliptic}, \cite{Xue_Zarhin_centers_Hodge_superelliptic}, 
\cite{Xue_Endomorphisms_Hodge_superelliptic}), their rational points (cf. \cite{Sutherland_counting_superelliptic}, \cite{Beneish_locally_soluble}, \cite{Legrand_twists_superelliptic}), their Galois representations (see e.g. \cite{Landesman_hyperelliptic_maximal_Galois},
\cite{Dokchitser_Arithmetic_hyperelliptic}, \cite{Anni_Dokchitser_Constructing_hyperelliptic} for hyperelliptic curves and \cite{goodman_superelliptic}, \cite{Pacetti_Villanueva_Galois_reps_superelliptic}, \cite{Bouw_Wewers_Computing_Lfcts_superelliptic} for superelliptic curves) and the isogenies between them (cf. \cite{Zarhin_non_isogenous_superelliptic}). 
Usually when considering the torsion or Galois representations of superelliptic Jacobians one excludes the primes dividing the exponent of $C$.
In this article we would like to study this ``exponential'' torsion.

To formulate the problem precisely, let $g := \dim J$. For any prime $\ell$ write $\rho_{J, \ell} : \Gal_K \to \Gl_{2g}(\ZZ_{\ell})$ for the Galois representation on $T_{\ell} J := \varprojlim J[\ell^n]$, the $\ell$-adic Tate module of~$J$. We are interested in the following question.
\begin{question} \label{q:galois_gps}
	What can be said about the image of $\rho_{J, \ell}$ for primes $\ell | m$?
\end{question}
Recall that a piece of $J[m]$ can be described using the roots of $f$.
Namely, there exists a subgroup $\Delta \subset J[m]$ such that $K(\Delta)$ is the splitting field of~$f$,
see \cite[Proposition~1]{Wawrow_torsion_superelliptic}. In the simplest case,
when $\GCD(m, r) = 1$, one has:
\begin{equation} \label{eqn:root_torsion}
	\Delta = \Span_{\ZZ}((\alpha_i, 0) - \infty : i =1, \ldots, r) \subset
	\Pic^0(C)[m],
\end{equation}
where $\alpha_1, \ldots, \alpha_r$ are the roots of $f$ and $\infty$ is unique point at infinity of $C$, see \cite[Proposition~1]{Wawrow_torsion_superelliptic}.
In case when $\ell = 2$, one has $\Delta = J[2]$. In particular, $\Gal(K(J[2])/K) = \Gal(f) \le S_r$. This observation allows one to answer Question~\ref{q:galois_gps} for a generic hyperelliptic curve.
For any ring epimorphism $R \to R/I$
and any subgroup $G \le \Gl_n(R)$, write
$G^{\mc S}$ for the preimage of a subset $\mc S \subset \Gl_n(R/I)$ in $G$. Note that $S_r$ may be treated as a subgroup of $\GSp_{2g}(\ZZ/2)$.
The above reasoning together with the existence of the Weil pairing easily implies that
\begin{equation} \label{eqn:image_2_adic}
	\im \rho_{J, 2} \subset \GSp_{2g}(\ZZ_2)^{\Gal(f)}.
\end{equation}
In case when $2|r$, $r \ge 6$ and $\Gal(f) = S_r$, this was shown to be an equality, cf. \cite[Corollary~1.2]{Yelton_Lifting_images}.\\

The goal of this article is to generalize~\eqref{eqn:image_2_adic} to the case when $m = \ell$ is an odd prime not dividing $r$. In the sequel, we consider only superelliptic Jacobians defined over $\QQ(\zeta_{\ell})$ for simplicity, but similar methods should apply for any number field containing $\zeta_{\ell}$. In other words, we use the following setup.
\begin{Setup} \label{setup}
	Let $\ell$ be an odd prime, $K := \QQ(\zeta_{\ell})$, $\mc O := \ZZ[\zeta_{\ell}]$ and $\lambda := 1- \zeta_{\ell}$. Denote by $J$ the Jacobian of the superelliptic
	curve $C$ defined by~\eqref{eqn:superelliptic}, where $m =\ell$ and $f \in K[x]$ is a monic separable polynomial of degree $r \ge 4$, $\ell \nmid r$. Also, denote the roots of $f$ by $\alpha_1, \ldots, \alpha_r$.
\end{Setup}
In this case, the set $\Delta$ is a proper subset of $J[\ell]$. It seems that determining the whole $\ell$-torsion explicitly is an impossible task.
There are however several known results for particular curves, see e.g. \cite{Tzermias_Explicit_rational}, \cite{Tzermias_cyclotomic_torsion}, \cite{Jedrzejak_torsion_superelliptic}, \cite{Jedrzejak_note_on_torsion}, \cite{arul_torsion_fermat_quotients}, \cite{Wawrow_torsion_superelliptic}.
The ``exponential'' torsion was also studied in \cite{Rasmussen_Tamagawa_cyclic_covers_ihara}, but in a different context.
A result of Arul from \cite{Arul_division} allows in principle to explicitly determine the $\ell$-torsion, but in practice the provided formulas are not well-suited for computations.
However, we can utilize the fact that $J$ has an additional structure of multiplication by $\mc O$. By comparing dimensions one checks that $\Delta = J[\lambda]$, see Section~\ref{sec:preliminaries}.
Moreover, $T_{\ell} J = \varprojlim J[\lambda^n]$ becomes a module over $\mc O_{\lambda} := \varprojlim \mc O/\lambda^n = \ZZ_{\ell}[\zeta_{\ell}]$. These observations will be crucial in the sequel.\\		

In order to describe the image of $\rho_{J, \ell}$, we need some notation concerning matrix groups. For any ring $R$, subgroup $G \le \Gl_n(R)$, and $X \subset R^{\times}$ we use the following notation:
\[ G_{\det \in X} := \{ A \in G : \det A \in X \}.  \]
Also, we use the notation $G^{\mc S}_{\det \in X}$ to denote the intersection $G^{\mc S} \cap G_{\det \in X}$.
Note that $\mc O_{\lambda}$ has a unique involution $z \mapsto \ol z$. Let $V$ be a free $\mc O_{\lambda}$-module
of finite rank and suppose that $\bb e : V \times V \to \mc O_{\lambda}$ is a non-degenerate Hermitian $\mc O_{\lambda}$-linear form. Denote by $\GU(V, \bb e)$ the associated general unitary group, i.e.
\[
\GU(V, \bb e) = \{ A \in \Gl(V) : \bb e(AP, AQ) = \mu \cdot \bb e(P, Q) \textrm{ for some } \mu \in \mc O_{\lambda}^{\times} \, \textrm{ and all } P, Q \in V \}.
\]
Recall that up to an isomorphism, there exist two non-degenerate Hermitian forms on~$V := \mc O_{\lambda}^d$, namely
$\bb e^+$ (given by the $d \times d$ identity matrix $I_d$) and $\bb e^-$ (given by the matrix $\diag(1, \ldots, 1, \alpha)$, where $\alpha \in \ZZ_{\ell}^{\times}$ is an arbitrary
non-square). We abbreviate $\GU_d^+(\mc O_{\lambda}) := \GU(\mc O_{\lambda}^d, \bb e^+)$ and $\GU_d^-(\mc O_{\lambda}) := \GU(\mc O_{\lambda}^d, \bb e^-)$.
In fact, if $2|d$, $\GU_d^+(\mc O_{\lambda}) \cong \GU_d^-(\mc O_{\lambda})$.
See Section~\ref{sec:preliminaries} for the details. Finally, denote by $\mc D_J$
the image of the map
\begin{equation*}
	\Omega_{J, \lambda} := \det {}_{\mc O_{\lambda}} \circ \rho_{J, \ell} : \Gal_K \to \mc O_{\lambda}^{\times},
\end{equation*}
where $\det {}_{\mc O_{\lambda}} : \Gl_{r-1}(\mc O_{\lambda}) \to \mc O_{\lambda}^{\times}$ denotes the $\mc O_{\lambda}$-linear determinant.
\begin{theorem} \label{thm:image_of_rho}
	Keep Setup~\ref{setup}. Then:
	\begin{equation} \label{eqn:main_inclusion}
		\rho_{J, \ell}(\Gal_K) \subset \GU_{r-1}^{\varepsilon}(\mc O_{\lambda})_{\det \in \mc D_J}^{\Gal(f)},
	\end{equation}
	where $\Gal(f) \subset S_r$ is considered as a subgroup of $\GU_{r-1}(\FF_{\ell})$ and
	$\varepsilon$ is the sign of the Legendre symbol $\left({\frac {r}{\ell}}\right)$. If moreover the following conditions hold:
	\begin{enumerate}[(1)]
		\item $\Gal(f)$ is a $2$-transitive subgroup of $S_r$,
		\item there exists a prime $\mf p$ of $\mc O$, $\mf p \nmid \ell$, such that $\ord_{\mf p}(\disc(f)) = 1$ and that the coefficients of $f$ are $\mf p$-integral,
	\end{enumerate}
	then the inclusion~\eqref{eqn:main_inclusion} is an equality.
\end{theorem}
Note that both conditions~(1) and~(2) hold for a generic polynomial $f \in K[x]$. Indeed,
the Galois group of a generic polynomial of degree $r$ is just the symmetric group~$S_r$.
Moreover, \cite[Corollary~28]{van_bommel_2014} proves that~(2) holds for most polynomials
in the case when $K = \QQ$ (the proof for number fields is similar). The idea behind the proof of Theorem~\ref{thm:image_of_rho} is to use that the reduction of $\rho_{J, \ell}(\Gal_K)$
modulo $\lambda$ equals $\Gal(f)$ and to observe that under some additional conditions this is the only subgroup of $\GU_{r-1}^{\varepsilon}(\mc O_{\lambda})_{\det \in \mc D_J}^{\Gal(f)}$ with this property. In the course of the proof we use also some additional information on $K(J[\lambda^2])$ coming from the descent theory. \\

Theorem~\ref{thm:image_of_rho} brings the study of the image of $\rho_{J, \ell}$ to studying $\mc D_J$.
One easily checks that $\mc D_J$ is contained~in
\[
\mc U_{\ell, r} := \{ d \in \mc O_{\lambda}^{\times} : d \cdot \ol{d} \in (1 + \ell \cdot (r-1) \cdot \ZZ_{\ell}) \}.
\]
Up to torsion, one may treat $\mc U_{\ell, r}$ as a $\ZZ_{\ell}$-lattice of rank $\frac{\ell + 1}{2}$, see Section~\ref{sec:endo_char}. 
It is natural to ask, whether the index of~$\mc D_J$ in $\mc U_{\ell, r}$ is finite. The following result
yields a positive answer.
\begin{theorem} \label{thm:DC}
	Keep the Setup~\ref{setup}. Then $[\mc U_{\ell, r} : \mu_{2 \ell} \mc D_J] = \ell^{\kappa}$ for some integer
	\[
	\kappa \le \ord_{\ell}(h_{\ell}^- \cdot c_{\ell, r}) - 1,
	\]
	where $h_{\ell}^-$ is the relative class number of $\QQ(\zeta_{\ell})$, $r_{\ell}$ is the multiplicative order of $r$ modulo~$\ell$ and
	\begin{equation} \label{eqn:c_lr_def}
		c_{\ell, r} := 
		\begin{cases}
			(r^{r_{\ell}} - 1)^{(\ell - 1)/2r_{\ell}}, & \textrm{ if } 2 \nmid r_{\ell},\\
			(r^{r_{\ell}/2} + 1)^{(\ell - 1)/r_{\ell}}, & \textrm{ if } 2 | r_{\ell}.
		\end{cases}
	\end{equation}
\end{theorem}
The strategy behind the proof of Theorem~\ref{thm:DC} is to relate $\Omega_{J, \lambda}$ to
a certain Hecke character with a specific infinity type. This allows one to construct an explicit $\ZZ_{\ell}$-lattice $\mc T_{\ell, r}$ inside of $\mc D_J$ and the question comes down to computing the determinant measuring the index of $\mc T_{\ell, r}$ in $\mc D_J$. It turns out that the determinant in question was computed by Hirabayashi in~\cite{Hirabayashi_generalization_Maillet_Demyanenko}. The matter of finding the torsion part of $\mc D_J$ comes down to studying $K(J[\lambda^2])$. In particular, if $J$ satisfies the assumptions~(1) and~(2) of Theorem~\ref{thm:image_of_rho}, then $\mu_{2 \ell} \subset \mc D_J$ (see Proposition~\ref{prop:torsion_in_DJ}).\\

The behaviour of the image of the Galois representation associated with an abelian variety is in general predicted by the Mumford--Tate conjecture. To every abelian variety~$A$ over a number field~$K$ one can associate a linear group $\MT(A)$
defined over $\QQ$, using the Hodge structure on the homology of~$A$. The Mumford--Tate conjecture states that for any prime $\ell$, 
the connected component of the Zariski closure of $\rho_{A, \ell}(\Gal_K)$ should be $\MT(A) \otimes \QQ_{\ell}$.
This conjecture serves as a link between two other significant problems -- the Hodge conjecture and the Tate conjecture.
The Mumford--Tate conjecture has been proven for many abelian varieties of types I, II and~III in Albert classification (cf. e.g. \cite{BGK_type_I}, \cite{BGK_type_III}, \cite{Farfan_torsion_type_iii}, \cite{gajda_hindry_remarks}; see also \cite{Farfan_Survey_HTMTC} for a survey).
The superelliptic Jacobians are of type IV in Albert classification. For this type fewer results are known (see e.g. \cite{Chi_IV} for a result concerning abelian varieties of prime dimension or \cite{Banaszak_Kaim_Garnek_Tate_module_type_IV} for some partial results). Also, it is known that superelliptic Jacobians of exponent~$3$,
whose defining polynomial has a large Galois group, satisfy the Mumford--Tate conjecture (cf. \cite[p. 103]{Zarhin_cyclic_covers} -- note that for abelian varieties the Tate conjecture implies the Mumford--Tate conjecture). The Hodge conjecture for superelliptic Jacobians was known in some particular cases, e.g. when $2 \ell | r$ (cf. \cite{Vasiu_some_cases_MT}) or if $r \not \equiv 1 \pmod{\ell}$ and $A_r \subset \Gal(f)$ (see \cite[Corollary 0.3]{Xue_Zarhin_Hodge_groups}).

It turns out that the Hodge, Tate and Mumford--Tate conjectures for the Jacobians in question are easy consequences
of the above results.
\begin{corollary} \label{cor:mt_conjecture}
	Keep the Setup~\ref{setup}. If $J$ satisfies the assumptions~(1) and~(2) of Theorem~\ref{thm:image_of_rho}, then
	$\End(J) = \mc O$. Moreover, the Hodge, Tate and Mumford--Tate conjectures hold for~$J$ and all its powers.
\end{corollary}
\noindent To our best knowledge, this is a new result. The proof of Corollary~\ref{cor:mt_conjecture} is given in Section~\ref{sec:mt_conjecture}. We end the introduction with a numerical example.
\begin{example} \label{ex:intro}
	Consider the superelliptic curve
	\[
	C : y^{11} = f = x^8 + x + 1
	\]
	over the field $K = \QQ(\zeta_{11})$. Then $J$ satisfies the assumptions~(1) and~(2) of Theorem~\ref{thm:DC} with $\mf p := 14731 \mc O$. In particular, the Hodge, Tate and Mumford--Tate conjectures hold for $J$. Moreover, $\left({\frac {8}{11}}\right) = -1$ and $\ord_{11}(h_{11}^- \cdot c_{11, 8}) = 1$. Using Proposition~\ref{prop:torsion_in_DJ} one checks also that $\mu_{22} \subset \mc D_J$.
	Thus $\mc D_J = \mc U_{11, 8}$ and
	\[
	\rho_{J, \ell}(\Gal_K) = \GU_7^-(\mc O_{\lambda})_{\det \in \mc U_{11, 8}}^{S_8}.
	\]
	In particular (see Example~\ref{ex:degree_division_field}):
	\begin{align*}
		[K(J[11]) : K] = 8! \cdot 11^{260}.
	\end{align*}
\end{example}
\subsection*{Outline of the paper} In Section~\ref{sec:preliminaries} we discuss
the needed preliminaries, including the superelliptic Jacobians and Hermitian pairings
over~$\mc O_{\lambda}$. Section~\ref{sec:lifting} is devoted to a proof
of a lifting result (cf. Proposition~\ref{prop:lifting}). This result
allows us to prove Theorem~\ref{thm:image_of_rho} in Section~\ref{sec:proof_of_thm_rho}.
Section~\ref{sec:endo_char} concerns the group $\mc D_J$, in particular
we prove Theorem~\ref{thm:DC}. In the last section we give a proof of
the Hodge, Tate and Mumford--Tate conjectures for superelliptic Jacobians satisfying the assumptions of Theorem~\ref{thm:image_of_rho}.
\subsection*{Acknowledgements}
The idea of proving Theorem~\ref{thm:image_of_rho} by using lifting methods was suggested by Davide Lombardo. The author would like to express gratitude for this concept. Additionally, the author extends thanks to Wojciech Gajda for many stimulating conversations. This project was started during the author's stay in Paris. The author wishes to express gratitude to his host, Marc Hindry, for the warm and welcoming hospitality extended during his stay. The author was supported by the research grant SONATINA 6 "The de Rham cohomology of $p$-group covers" UMO-2022/44/C/ST1/00033
awarded by National Science Centre, Poland. 
	
\section{Preliminaries} \label{sec:preliminaries}
The goal of this section is to introduce the needed notation and to discuss some preliminaries on the superelliptic Jacobians.
In particular, we prove the inclusion~\eqref{eqn:main_inclusion}. For any number field $L$ we write $\mc O_L$ for its ring
of integers and $\disc(L/\QQ)$ for its discriminant over $\QQ$. For any prime ideal $\mf p$ of $\mc O_L$ we denote:
\begin{itemize}
	\item the associated valuation by $\ord_{\mf p}$,
	\item the completion of $L$ with respect to $\mf p$ by $L_{\mf p}$,
	\item the valuation ring of $L_{\mf p}$ by $\mc O_{L, \mf p}$,
	\item an arbitrary uniformizer of $\mf p$ by $\pi_{\mf p}$.
\end{itemize}
Similarly, if $v$ is a place of $L$, we write $L_v$ for the corresponding completion.
If $v$ is finite, we use the notation~$\mf p_v$ for the corresponding
prime ideal of $\mc O_L$ and put $\mc O_{L, v} := \mc O_{L, \mf p_v}$.
Also, let $\mc O_{L, v} := L_v$ for any infinite place $v$.
For any polynomial~$f$ we write $\disc(f)$ for its discriminant.\\

From now on, we use the Setup~\ref{setup}. In this case we abbreviate $\mc O_{K, v}$ to $\mc O_v$, etc. Note that $\lambda$ is a prime element of $\mc O$. We identify it with the prime ideal $\lambda \mc O$. Moreover, $\ell \mc O = \lambda^{\ell - 1} \mc O$.
In the sequel, we use several times the congruence:
\begin{equation} \label{eqn:ol_lambda_is_minus_lambda}
	\ol{\lambda} \equiv - \lambda \pmod{\lambda^2},
\end{equation}
which can be proven as follows:
\begin{align*}
	\ol{\lambda} &= 1 - \zeta_{\ell}^{\ell - 1} = \frac{1 - \zeta_{\ell}^{\ell - 1}}{1 - \zeta_{\ell}}
	\cdot \lambda\\
	&=  (1 + \zeta_{\ell} + \ldots + \zeta_{\ell}^{\ell - 2})
	\cdot \lambda \equiv (\ell - 1) \cdot \lambda \equiv - \lambda \pmod{\lambda^2}.
\end{align*}
For a matrix $\bb A \in \Gl_d(\mc O_{\lambda})$ we often write
$\det_{\mc O_{\lambda}} \bb A$ to denote the determinant of $\bb A$ considered as an $\mc O_{\lambda}$-linear map
and $\det_{\ZZ_{\ell}} \bb A$ when we consider $\bb A$ to be an element of $\Gl_{d \cdot (\ell - 1)}(\ZZ_{\ell})$.
Note that $\det_{\ZZ_{\ell}} \bb A = N_{K_{\lambda}/\QQ_{\ell}}(\det_{\mc O_{\lambda}} \bb A)$.\\

\noindent Recall that the map
\[
[\zeta_{\ell}] : C \to C, \quad (x, y) \mapsto (x, \zeta_{\ell} \cdot y),
\]
is an automorphism of order $\ell$. This leads to a map $\mc O \to \End(J)$, which turns out to be an embedding (see e.g. \cite[Section~1]{Zarhin_endomorphisms_superelliptic}). In what follows, we use often the fact that the Rosati involution on $\mc O$ is given by the complex conjugation. The dimension of $J$, denoted $g$, equals $\frac{1}{2}(r - 1) \cdot (\ell - 1)$ (this follows e.g. from \cite[Proposition~3.7.3 (c)]{Stichtenoth_algebraic_function_fields}). Therefore, by~\cite[Proposition~2.2.1]{Ribet_galois_action_real_multiplication},
for every $n \ge 1$:
\[
J[\lambda^n] \cong (\mc O/\lambda^n)^{2 \dim J/(\ell - 1)} \cong (\mc O/\lambda^n)^{r - 1}.
\]
In particular, $T_{\ell} J \cong \mc O_{\lambda}^{r - 1}$. Observe that since $[\zeta_{\ell}] (\alpha_i, 0) = (\alpha_i, 0)$, the group $\Delta$ defined by~\eqref{eqn:root_torsion} is an $\FF_{\ell}$-linear subspace of $J[\ell]$.
Moreover, $\dim_{\FF_{\ell}} \Delta = r - 1$,
since $\Delta$ is spanned by the divisors $(\alpha_i, 0) - \infty$ for $i = 1, \ldots, r$ with the only relation:
\[
\sum_{i = 1}^r ((\alpha_i, 0) - \infty) = \divv(y) \sim 0 \quad \textrm{ in } \Pic^0(C),
\]
cf. \cite[Proposition~1]{Wawrow_torsion_superelliptic}. Hence $J[\lambda] = \Delta$. In particular,
$\Gal(f) = \Gal(K(J[\lambda])/K) \le \Gl_{r-1}(\FF_{\ell})$. Note that this inclusion
is compatible with $\Gal(f) \le S_r$ and $S_r \le \Gl_{r-1}(\FF_{\ell})$,
where the second inclusion comes from the permutation action on $\FF_{\ell}^r/\diag(\FF_{\ell})$.

We study now the Weil pairing on $T_{\ell} J$. To this end we quote the general result on Weil pairings on abelian varieties
with multiplications by an order in some number field. We give also the proof based on \cite{Zarhin_on_MO} for a lack of better reference. For any Galois module~$M$, write
$M(1)$ for its twist by the cyclotomic character.
\begin{lemma} \label{lem:weil_general}
	Let $A$ be a principally polarized abelian variety over a number field~$L$.
	Suppose that a number field $M$ embedds into $\End(A) \otimes \QQ$. Let $\OO := M \cap \End(A)$.
	Let also $\OO^*$ denote the different ideal of $\OO$ over $\ZZ$. For any natural number $n$ there exists a unique skew-Hermitian (with respect to the Rosati involution), $\OO$-linear on the second coordinate,
	$\Gal_L$-equivariant and non-degenerate pairing
	\[
	\mf e_n : A[n] \times A[n] \to (\OO^*/n)(1)
	\]
	(where $\OO^*/n$ is considered with the trivial Galois action) such that:
	\begin{itemize}
		\item $\tr_{M/\QQ} \mf e_n$ equals the standard Weil pairing:
		\[
		e_n : A[n] \times A[n] \to (\ZZ/n)(1),
		\]
		
		\item for any $n_1, n_2$, $P \in A[n_1 \cdot n_2]$ and
		$Q \in A[n_1]$ one has:
		\[
		\mf e_{n_1 n_2}(P, Q) = \mf e_{n_1}([n_2] P, Q).
		\]
	\end{itemize}
	Similarly, for any prime $\ell$ there exists a unique skew-Hermitian (with respect to the Rosati involution), $(\OO^* \otimes \ZZ_{\ell})$-linear on the second variable
	and $\Gal_L$-equivariant non-degenerate pairing
	\[
	\mf e_{\ell^{\infty}} : T_{\ell} A \times T_{\ell} A \to (\OO^* \otimes \ZZ_{\ell})(1)
	\]
	(where $\OO^* \otimes \ZZ_{\ell}$ is considered with the trivial Galois action) such that $\tr_{M/\QQ} \mf e_{\ell^{\infty}}$ equals the standard Weil pairing
	$e_{\ell^{\infty}} : T_{\ell} A \times T_{\ell} A \to \ZZ_{\ell}(1)$.
\end{lemma}
\begin{proof}
	In the proof we will often use the fact that the map
	\begin{equation} \label{eqn:different_and_isomorphism_of_homs}
		\OO^*/n \to \Hom_{\ZZ}(\OO/n, \ZZ/n), \qquad a \mapsto f_a, 
	\end{equation}
	where $f_a(x) := \tr_{M/\QQ}(a \cdot x)$, is an isomorphism of $\ZZ$-modules.
	Indeed, by definition $\OO^* \cong \Hom_{\ZZ}(\OO, \ZZ)$
	(see e.g. \cite[Definition~III.2.1]{Neukirch_ANT}).
	Moreover, since $\OO$ and $\ZZ$ are finitely generated projective $\ZZ$-modules,
	$\Hom_{\ZZ}(\OO, \ZZ) \otimes \ZZ/n \cong \Hom_{\ZZ}(\OO/n, \ZZ/n)$.\\
	
	Denote the Rosati involution on $\OO$ by $z \mapsto \ol z$.
	Fix $P, Q \in A[n]$ and consider the map:
	\[
	\OO/n \to \ZZ/n, \qquad a \mapsto e_n(P, a \cdot Q).
	\]
	Since~\eqref{eqn:different_and_isomorphism_of_homs} is an isomorphism, there exists a unique element $\mf e_n(P, Q) \in \OO^*/n$ such that:
	\[
	e_n(P, a \cdot Q) = \tr_{M/\QQ}(\mf e_n(P, Q) \cdot a) \qquad \textrm{ for any } a \in \OO/n.
	\]
	The map $(P, Q) \mapsto \mf e_n(P, Q)$ is $\ZZ$-bilinear and $\OO$-linear in the second coordinate. Indeed, for example we have:
	\begin{align*}
		\tr_{M/\QQ}((\mf e_n(P_1, Q) + \mf e_n(P_2, Q)) \cdot a)
		&= \tr_{M/\QQ}(\mf e_n(P_1, Q) \cdot a) + \tr_{M/\QQ}(\mf e_n(P_2, Q) \cdot a)\\
		&= e_n(P_1, a \cdot Q) + e_n(P_2, a \cdot Q) = e_n(P_1 + P_2, a \cdot Q)\\
		&= \tr_{M/\QQ}(\mf e_n(P_1 + P_2, Q) \cdot a),
	\end{align*}
	which implies that $\mf e_n(P_1 + P_2, Q) = \mf e_n(P_1, Q) + \mf e_n(P_2, Q)$.
	We prove now that $\mf e_n$ is skew-Hermitian. Indeed, this follows from the fact that $a$ and $\ol a$
	are dual maps with respect to $e_n$:
	\begin{align*}
		\tr_{M/\QQ}(\mf e_n(P, Q) \cdot a) &= e_n(P, a \cdot Q)
		= - e_n(a \cdot Q, P) = - e_n(Q, \ol a \cdot P)\\
		&= - \tr_{M/\QQ}(\mf e_n(Q, P) \cdot \ol a) = - \tr_{M/\QQ}(\ol{\mf e_n(Q, P)} \cdot a).
	\end{align*}
	But since~\eqref{eqn:different_and_isomorphism_of_homs} is an isomorphism, we obtain that $\mf e_n(P, Q) = - \ol{\mf e_n(Q, P)}$. Similarly one shows the remaining properties of $\mf e_n$.
	The pairing $\mf e_{\ell^{\infty}}$ is constructed as the inverse limit of
	pairings $\mf e_{\ell^n}$.
\end{proof}
We switch now back to the notation of the article.
\begin{corollary} \label{cor:weil_pairing_superelliptic}
	Keep Setup~\ref{setup}. There exists a Hermitian, $\mc O_{\lambda}$-linear on the second coordinate, non-degenerate pairing of $\Gal_K$-modules:
	\[
	\bb e_{\lambda^{\infty}} : T_{\ell} J \times T_{\ell} J \to \mc O_{\lambda}(1).
	\]
\end{corollary}
\begin{proof}
	Using Lemma~\ref{lem:weil_general} for $L = M = K$, we obtain the $\mc O$-linear Weil pairing
	\[
	\mf e_{\lambda^{\infty}} : T_{\ell} J \times T_{\ell} J \to (\mc O^* \otimes \ZZ_{\ell}) (1).
	\]
	Recall that $\mc O^* = \lambda^{-(\ell - 2)} \cdot \mc O$ (see e.g. \cite[exercise~III.\S 2.5]{Neukirch_ANT}). Define:
	\[
	\bb e_{\lambda^{\infty}} := (\lambda^{\ell - 2} - \ol{\lambda}^{\ell - 2}) \cdot \mf e_{\lambda^{\infty}}.
	\]
	Then $\bb e_{\lambda^{\infty}}(P, Q) \in \mc O_{\lambda}$ for any $P, Q \in T_{\ell} J$, since $\ord_{\lambda}(\lambda^{\ell - 2} - \ol{\lambda}^{\ell - 2}) = \ell - 2$ by~\eqref{eqn:ol_lambda_is_minus_lambda}.
	Moreover,
	since $\mf e_{\ell}$ is skew-Hermitian and $\lambda^{\ell - 2} - \ol{\lambda}^{\ell - 2}$
	is an imaginary number, $\bb e_{\lambda^{\infty}}$ is a Hermitian pairing.
\end{proof}
We discuss now briefly the theory of Hermitian forms over $\mc O_{\lambda}$. Suppose that $V$ is a
free $\mc O_{\lambda}$-module of rank~$d$ and that
\begin{equation*}
	\bb e : V \times V \to \mc O_{\lambda}
\end{equation*}
is a non-degenerate $\mc O_{\lambda}$-linear Hermitian form. In this case we define $\GU(V, \bb e)$ as in Section~\ref{sec:intro}. Note that the multiplier $\mu$ yields
a well-defined homomorphism:
\[
\mu : \GU(V, \bb e) \to \mc O_{\lambda}^{\times}.
\]
The kernel of $\mu$ is denoted by $\U(V, \bb e)$ and the kernel of $\det_{\mc O_{\lambda}} : \U(V, \bb e) \to \mc O_{\lambda}^{\times}$ by $\SU(V, \bb e)$.
Note that if $\Gamma$ is the matrix of $\bb e_{\lambda^{\infty}}$, then
for any $\bb A \in \GU(V, \bb e)$ by taking the $\mc O_{\lambda}$-linear determinant of the equation
\[
\bb A^{\dagger} \cdot \Gamma \cdot \bb A = \mu(A) \cdot \Gamma,
\]
(where $\bb A^{\dagger}$ denotes the conjugate transpose of a matrix $\bb A$) we obtain:
\begin{equation} \label{eqn:ol_det_det_equals_mu_r-1}
	\ol{\det {}_{\mc O_{\lambda}} \bb A} \cdot \det {}_{\mc O_{\lambda}} \bb A = \mu(\bb A)^d.
\end{equation}
Also, we write $\GU(V/\lambda^n, \bb e)$ for the group of automorphisms
of the free $\mc O/\lambda^n$-module $V/\lambda^n$ preserving up to a scalar the form $\bb e_{\lambda^n}$, the reduction of $\bb e$ modulo $\lambda^n$.
In an analogous way one defines
$\U(V/\lambda^n, \bb e)$ and $\SU(V/\lambda^n, \bb e)$.\\

It turns out that if $\bb e_{\lambda}$ is non-degenerate, it uniquely determines $\bb e$, cf. \cite[Theorem~2]{Hermitian_skew_hermitian_local_rings}.
The form $\bb e_{\lambda}$ is then a symmetric form (since the involution $z \mapsto \ol z$ is trivial on $\FF_{\ell}$).
Recall that the determinant of a pairing is defined as the determinant of its matrix in any basis.
It is a well-defined invertible element up to squares in the given ring.
In case of finite fields, the determinant is the complete invariant for non-degenerate symmetric forms of the same dimension (cf. \cite[\S 1.7, Proposition~5]{Serre_course_arithmetic}). 
Thus there exist exist two non-degenerate symmetric bilinear forms on $\FF_{\ell}^d$, namely
$\bb e^+_{\lambda}$ (given by the identity matrix $I_d$) and $\bb e^-_{\lambda}$ (given by the matrix $\diag(1, \ldots, 1, \alpha)$, where $\alpha \in (\FF_{\ell})^{\times}$ is an arbitrary
non-square residue). By the above discussion, the non-degenerate bilinear Hermitian forms on $\mc O_{\lambda}^d$ are classified by the class of the determinant 
in $\mc O_{\lambda}^{\times}/(\mc O_{\lambda}^{\times})^2$. Therefore there are two isomorphism classes of such forms, represented by
the forms $\bb e^+$ and $\bb e^-$ with matrices $I_d$ and $\diag(1, \ldots, 1, \alpha)$ (where $\alpha \in \mc O_{\lambda}^{\times}$ is an arbitrary non-square) respectively.
We abbreviate $\GU^+_d(\mc O_{\lambda}) := \GU(\mc O_{\lambda}^d, \bb e^+)$, $\GU^-_d(\mc O_{\lambda}) := \GU(\mc O_{\lambda}^d, \bb e^-)$, etc.
In case when $2 \nmid d$, the Hermitian form given by the matrix $\alpha \cdot I_d$ has a non-square determinant, thus it is equivalent
to $\bb e^-$. Since an automorphism $\bb A \in \Gl_d(\mc O_{\lambda})$ commutes with the form given by the matrix $\alpha \cdot I_d$ if and only if it commutes
with the form given by $I_d$, we obtain $\GU^+_d(\mc O_{\lambda}) \cong \GU^-_d(\mc O_{\lambda})$.
\begin{lemma}
	We have $\rho_{J, \ell}(\Gal_K) \subset \GU_{r-1}^{\varepsilon}(\mc O_{\lambda})$, where $\varepsilon$ is the sign of the Legendre symbol $ \left({\frac {r}{\ell}}\right)$.
\end{lemma}
\begin{proof}
	Clearly, $\rho_{J, \ell}(\Gal_K) \subset \GU(T_{\ell} J, \bb e_{\lambda^{\infty}})$. Thus
	by the previous discussion, the proof is immediate if $2 | r$. Therefore, from now on we suppose that $r$ is odd.
	
	We start by proving the result in the special case when $\Gal(f)$ is a $2$-transitive subgroup of~$S_r$. Under this assumption we will show that $\bb e_{\lambda^{\infty}}$ is equivalent to $\bb e^{\varepsilon}$. Let $P_i := (\alpha_i, 0) - \infty \in J[\lambda]$ for $i = 1, \ldots, r$.
	Denote $c := \bb e_{\lambda}(P_1, P_2)$. By the assumption, for any $i \neq j$ there exists $\sigma \in \Gal(f)$ such that $\sigma(P_1) = P_i$, $\sigma(P_2) = P_j$.
	Thus, since $\bb e_{\lambda}$ is Galois equivariant and $\zeta_{\ell} \in K$:
	\[
	\bb e_{\lambda}(P_i, P_j) = \bb e_{\lambda}(\sigma(P_1), \sigma(P_2)) = \bb e_{\lambda}(P_1, P_2) = c.
	\]
	Moreover, for any $1 \le i \le r$, since $P_1 + \ldots + P_r = 0$:
	\begin{align*}
		\bb e_{\lambda}(P_i, P_i) = \bb e_{\lambda} \left(P_i, - \sum_{\substack{j = 1 \\ j \neq i}}^r P_j \right) = -(r-1) \cdot c.
	\end{align*}
	Therefore, in the basis $P_1, \ldots, P_{r-1}$ the matrix of $\bb e_{\lambda}$ is given by $c \cdot (E - r \cdot I_{r-1})$,
	where the matrix $E \in M_{r-1}(\FF_{\ell})$ is entirely composed of ones. The matrix $E$ is of rank $1$ and its eigenvalues are $r-1, 0, 0, \ldots, 0$.
	Thus the eigenvalues of $c \cdot (E - r \cdot I_{r-1})$ are $-c, -c \cdot r, \ldots, - c \cdot r$ and $\det(c \cdot (E - r \cdot I_{r-1}))
	= c^{r-1} \cdot r^{r-2}$. In particular, the determinant of $\bb e_{\lambda}$ is the class of $r$ in $\FF_{\ell}^{\times}/(\FF_{\ell}^{\times})^2$. This ends the proof in this case.\\
	
	In general, if $\Gal(f)$ is an arbitrary subgroup of~$S_r$, let $\mc K := K(a_0, \ldots, a_r)$ be the function field in $r+1$ variables.
	Write $\mc O_{\mc K} := \QQ[a_0, \ldots, a_r]$. Let $\bb f_{\bb a}(x) := a_0 + a_1 x + \ldots + a_r x^r \in \mc K[x]$.
	Consider the family of curves $\bb C \to \Spec(\mc O_{\mc K})$ given by $y^{\ell} = \bb f(x)$.
	Let $\mc J \to \Spec(\mc O_{\mc K})$ be the relative Jacobian. Write $J_{\eta}$ for the generic fiber of this family.
	Note that $\Gal(\bb f) \cong S_r$ (see e.g. \cite[Example VI.\S 2.4, p. 272]{Lang_Algebra}).
	Therefore, by the above discussion, $\rho_{J_{\eta}, \lambda}(G_{\mc K}) \subset \GU_{r-1}^{\varepsilon}(\mc O_{\lambda})$. On the other hand,
	since $\mc J$ specializes to $J$, we have:
	\[
	\rho_{J, \ell}(\Gal_K) \subset \rho_{J_{\eta}, \lambda}(\Gal_{\mc K})
	\]
	(see e.g. \cite[Lemma 6.1.6~(b)]{cantoral_lombardo_monodromy_QM_jacobians}). This ends the proof.
\end{proof}
Corollary~\ref{cor:weil_pairing_superelliptic} implies also that
\begin{equation*} 
	\mu \circ \rho_{J, \ell} = \chi_{cyc} : \Gal_K \to \ZZ_{\ell}^{\times}
\end{equation*}
is the cyclotomic character. In particular, the equation~\eqref{eqn:ol_det_det_equals_mu_r-1} becomes
\begin{equation} \label{eqn:ol_det_det_is_cyc}
	\ol{\det{}_{\mc O_{\lambda}} (\rho_{J, \ell}(\sigma))} \cdot \det{}_{\mc O_{\lambda}} (\rho_{J, \ell}(\sigma)) = \chi_{cyc}^{r - 1}(\sigma) \quad \textrm{ for every } \sigma \in \Gal_K.
\end{equation}
\section{A lifting result} \label{sec:lifting}
When computing images of abelian varieties, one often encounters
lifting problems of the following form.
\begin{question} \label{q:lifting_question}
	Let $R$ be a topological ring with a quotient $R/\mc I$. Assume that
	$\mc G$, $\mc H$ are closed subgroups of $\Gl_d(R)$ such that
	$\mc H \le \mc G$ and that the images of $\mc G$ and $\mc H$ in $\Gl_d(R/\mc I)$
	are equal. Does it follow that $\mc G = \mc H$?
\end{question}
The answer is positive under some additional assumptions if for instance $R$ is the ring of Witt vectors over a finite field of characteristic~$\ell$, $\mc I = \ell R$ and $\mc G$ is an algebraic group (see \cite{Vasiu_surjectivity_criteria_i} and \cite{Weigel_profinite_completion})
or if $R = \ZZ_{\ell}$, $\mc I = \ell R$ and $\mc G$ is an intersection of an algebraic group and preimage of a subgroup of $\Gl_d(\FF_{\ell})$ under the reduction map (cf. \cite{Jones_Gl2_reps} and \cite{Yelton_Lifting_images}).
In our case, we would like to take $R := \mc O_{\lambda}$ and $\mc G := \SU_d^{\varepsilon}(\mc O_{\lambda})^{\Gal(f)}$. Thus there arise the following two difficulties:
\begin{itemize}
	\item $\mc O_{\lambda}$ is a ramified extension of $\ZZ_{\ell}$,
	
	\item $\mc G$ is an intersection of a preimage of a subgroup of $\Gl_d(\FF_{\ell})$
	with a group $\SU_d^{\varepsilon}(\mc O_{\lambda})$, which is not algebraic over $\mc O_{\lambda}$.
\end{itemize}
However, it turns out that in this context the answer to Question~\ref{q:lifting_question}
is still positive if we take $\mc I := \lambda^2 \mc O_{\lambda}$. The following is
the main result of this section.
\begin{proposition} \label{prop:lifting}
	Suppose that $\bb e : V \times V \to \mc O_{\lambda}$ is a Hermitian
	form on a free $\mc O_{\lambda}$-module $V$ of rank $d \ge 3$. Let $\mc S$ be a subgroup of $\SU(V/\lambda, \bb e_{\lambda})$.
	Suppose that $\mc H \le \SU(V, \bb e)^{\mc S}$ is a closed subgroup such that its reduction modulo $\lambda^2$ is $\SU(V/\lambda^2, \bb e)^{\mc S}$.
	Then $\mc H = \SU(V, \bb e)^{\mc S}$.
\end{proposition}
\begin{remark}
	Observe that Proposition~\ref{prop:lifting} is not true with~$\lambda^2$ replaced by~$\lambda$. Indeed, take for example $\mc H := \SO_d(\ZZ_{\ell})$ (the special orthogonal group) and $\mc G := \SU_d^+(\mc O_{\lambda})$. Then both $\mc G$ and $\mc H$ surject onto 
	$\SU_d^+(\FF_{\ell})$, but clearly $\mc H \neq \mc G$.
\end{remark}
\begin{remark} \label{rmk:unitary_is_algebraic}
	Note that even though $\SU(V, \bb e)$ is not an algebraic group over $\mc O_{\lambda}$,
	it can be treated as an algebraic group over $\mc O_{\lambda}^+ := \ZZ_{\ell}[\zeta_{\ell} + \zeta_{\ell}^{-1}]$ via a ``Weil restriction''.
	In other words, it is easy to check that the group
	\[
	\{ (X, Y) \in M_d(\mc O_{\lambda}^+) \times M_d(\mc O_{\lambda}^+) : X + \zeta_{\ell} Y \in \SU_d(V, \bb e) \}
	\]
	is algebraic.
\end{remark}
In what follows,
matrices in $M_d(\mc O/\lambda^n)$ are written in bold face and
matrices in $M_d(\FF_{\ell})$ are written in regular font. For any matrix
$\bb A \in M_d(\mc O/\lambda^n)$, we denote by $A_0, \ldots, A_{n-1} \in M_d(\FF_{\ell})$ its ``digits'' in the ``$\lambda$-adic expansion'':
\begin{equation} \label{eqn:adic_exp_of_matrix}
	\bb A = A_0 + \lambda \cdot A_1 + \ldots + \lambda^{n-1} \cdot A_{n-1}.
\end{equation}
For any group $G \le \Gl_d(\mc O/\lambda^n)$ and $1 \le k \le n$, denote by $G_k$ the kernel of
the reduction map $G \to \Gl_d(\mc O/\lambda^k)$. We use a similar notation for groups $G \le \Gl_d(\mc O_{\lambda})$.
Note that $\bb A \in \Gl_d(\mc O/\lambda^n)_k$ for $1 \le k \le n$
if and only if $A_0 = I_d$ and $A_1 = \ldots = A_{k-1} = 0$. Moreover, if $\bb A = I_d + \lambda^{n-1} \cdot A_{n-1} \in \Gl_d(\mc O/\lambda^n)_{n-1}$,
then
\begin{equation} \label{eqn:trace_and_det}
	\det {}_{\mc O_{\lambda}}(\bb A) = 1 + \lambda^{n-1} \cdot \tr(A_{n-1}).
\end{equation}
Indeed, this follows easily from the fact that $\det(I_d - x A_{n-1})$ is the reciprocal of the characteristic polynomial
of $A_{n-1}$.

The rest of this section may be skipped at first reading.
We start to work on the proof of Proposition~\ref{prop:lifting}.
From the discussion in Section~\ref{sec:preliminaries} we may choose a basis
of $V$ in which~$\bb e$ has matrix of the form
$\Gamma = \diag(\alpha_1, \ldots, \alpha_d)$, where $\alpha_1, \ldots, \alpha_d \in \ZZ_{\ell}^{\times}$
(actually, one can assume also that $\alpha_1 = \ldots = \alpha_{d-1} = 1$, but we won't need this).
Thus we may identify
\[
\SU(V, \bb e) =  \{ \bb A \in \Sl_d(\mc O_{\lambda}) : \bb A^{\dagger} \Gamma \bb A = \Gamma \}.
\]

Write $\llbracket \cdot, \cdot \rrbracket$ for the commutator in the group $\Gl_d$ (i.e. $\llbracket A, B \rrbracket := A B A^{-1} B^{-1}$)
and $[\cdot, \cdot]$ for the Lie bracket in~$M_d$ (i.e. $[A, B] := AB - BA$).
\begin{lemma} \label{lem:computation_of_commutator}
	Let $N := \lfloor \frac{n-1}{2} \rfloor$. 
	Suppose that $\bb A, \bb B \in \Gl_d(\mc O/\lambda^n)_N$. Then:
	\[
	\llbracket \bb A, \bb B \rrbracket =
	\begin{cases}
		I_d + \lambda^{n-1} [A_N, B_N], & 2 \nmid n\\
		I_d + \lambda^{n-2} [A_N, B_N] + \lambda^{n-1} ([A_N, B_{N+1}] + [A_{N+1}, B_N]), & 2 | n, n \ge 6\\
		I_d + \lambda^2 [A_1, B_1] + \lambda^3 ([A_1, B_2] + [A_2, B_1] + [B_1, A_1] A_1 + [B_1, A_1] B_1), & n = 4.
	\end{cases}
	\]
\end{lemma}
\begin{proof}
	For any given $n$ one can verify this identity, using a simple Sage Math code, see the Appendix.
	In particular, this proves the result for $n = 4$. For the remaining $n$ we split the proof into two cases:
	\begin{itemize}[leftmargin=*]
		\item \bb{Case I:} $n$ is odd.
		\item[] In this case $n = 2N+1$. Note that for any $\bb X, \bb Y \in \Gl_d(\mc O/\lambda^n)_N$:
		\[
		\bb X \cdot \bb Y = I_d + \sum_{i = N}^{n-1} \lambda^i \cdot (X_i + Y_i) + \lambda^{n-1} \cdot X_N Y_N.
		\]
		The equation $\bb X \cdot \bb Y = I_d$ is equivalent to $(\bb X \cdot \bb Y)_i = 0$ for $i = N, \ldots, n-1$.
		Thus, in particular:
		\[
		\bb X^{-1} = I_d - \sum_{i = N}^{n-1} \lambda^i \cdot X_i + \lambda^{n-1} \cdot X_N^2.
		\]
		This yields:
		\begin{align*}
			\llbracket \bb A, \bb B \rrbracket &= \bb A \cdot \bb B \cdot \bb A^{-1} \cdot \bb B^{-1} = \left(I_d + \sum_{i = N}^{n-1} \lambda^i A_i \right) \cdot \left(I_d + \sum_{i = N}^{n-1} \lambda^i B_i \right)\\
			&\cdot \left(I_d - \sum_{i = N}^{n-1} \lambda^i \cdot A_i + \lambda^{n-1} \cdot A_N^2 \right) 
			\cdot \left(I_d - \sum_{i = N}^{n-1} \lambda^i \cdot B_i + \lambda^{n-1} \cdot B_N^2 \right)
		\end{align*}
		\begin{align*}
			&= \left(I_d + \sum_{i = N}^{n-1} \lambda^i \cdot (A_i + B_i) + \lambda^{n-1} \cdot A_N B_N \right)\\
			&\cdot \left(I_d - \sum_{i = N}^{n-1} \lambda^i \cdot (A_i + B_i) + \lambda^{n-1} \cdot (A_N^2 + B_N^2 + A_N B_N) \right)\\
			&= I_d + \lambda^{n-1} \cdot [A_N, B_N].
		\end{align*}
		\item \bb{Case II:} $n \ge 6$ is even.
		\item[] In this case $n = 2N+2$. Similarly as before, note that for any $\bb X, \bb Y \in \Gl_d(\mc O/\lambda^n)_N$:
		\begin{align*}
			\bb X \cdot \bb Y &= 
			I_d + \sum_{i = N}^{n-1} \lambda^i \cdot (X_i + Y_i) + \lambda^{n-2} \cdot X_N Y_N + \lambda^{n-1} \cdot (X_N Y_{N+1}
			+ X_{N+1} Y_N),\\
			\bb X^{-1} &=
			I_d - \sum_{i = N}^{n-1} \lambda^i \cdot X_i
			+ \lambda^{n-1} \cdot X_N^2 + \lambda^{n-1} \cdot
			(X_N X_{N+1} + X_{N+1} X_N)
		\end{align*}
		(we used here the fact that $N + 1 < n-2$, since $n \ge 6$). Thus, by similar computations:
		\begingroup
		\allowdisplaybreaks
		\begin{align*}
			\llbracket \bb A, \bb B \rrbracket &= \left(I_d + \sum_{i = N}^{n-1} \lambda^i A_i \right) \cdot \left(I_d + \sum_{i = N}^{n-1} \lambda^i B_i \right)\\
			&\cdot \left(I_d - \sum_{i = N}^{n-1} \lambda^i \cdot A_i + \lambda^{n-2} \cdot A_N^2 + \lambda^{n-1} \cdot
			(A_N A_{N+1} + A_{N+1} A_N) \right)\\ 
			&\cdot \left(I_d - \sum_{i = N}^{n-1} \lambda^i \cdot B_i + \lambda^{n-2} \cdot B_N^2 + \lambda^{n-1} \cdot
			(B_N B_{N+1} + B_{N+1} B_N)\right)\\
			&= \left(I_d + \sum_{i = N}^{n-1} \lambda^i \cdot (A_i + B_i) + \lambda^{n-2} A_N B_N + \lambda^{n-1} \cdot (A_N B_{N+1}
			+ A_{N+1} B_N)\right)\\
			&\cdot \bigg(I_d - \sum_{i = N}^{n-1} \lambda^i \cdot (A_i + B_i) + \lambda^{n-2} \cdot (A_N^2 + B_N^2 + A_N B_N)\\
			&+\lambda^{n-1} \cdot ((A_N A_{N+1} + A_{N+1} A_N) + (B_N B_{N+1} + B_{N+1} B_N) + A_N B_{N+1} + A_{N+1}B_N) \bigg)\\
			&=I_d  + \lambda^{n-2} \cdot (2 A_N B_N + A_N^2 + B_N^2 - (A_N + B_N)^2)\\
			&+ \lambda^{n-1} \cdot ((A_N B_{N+1}
			+ A_{N+1} B_N) + (A_N A_{N+1} + A_{N+1} A_N + B_N B_{N+1} + B_{N+1} B_N\\
			& + A_N B_{N+1} + A_{N+1}B_N) - (A_N + B_N)(A_{N+1} + B_{N+1}) - (A_{N+1} + B_{N+1})(A_N + B_N))\\
			&= I_d + \lambda^{n-2} \cdot [A_N, B_N] + \lambda^{n-1} \cdot ([A_N, B_{N+1}] + [A_{N+1}, B_N]).
		\end{align*}
		\endgroup
		This ends the proof. \qedhere
	\end{itemize}
\end{proof}
\begin{lemma} \label{lem:gl_in_center}
	For any $1 \le i \le j \le n$, $i + j \ge n$, the subgroup $\Gl_d(\mc O/\lambda^n)_j$ is in the center of $\Gl_d(\mc O/\lambda^n)_i$.
\end{lemma}
\begin{proof}
	If $\bb A \in \Gl_d(\mc O/\lambda^n)_i$, $\bb B \in \Gl_d(\mc O/\lambda^n)_j$,
	then (since $i + j \ge n$)
	\begin{align*}
		\bb A \cdot \bb B = I_d + \sum_{i = 1}^{n-1} \lambda^i \cdot (A_i + B_i) = \bb B \cdot \bb A.
	\end{align*}
	This ends the proof.
\end{proof}
Recall that if $G \le \Gl_d(\mc O_{\lambda})$ is an algebraic group with a Lie algebra $\Lie(G) \le M_d(\mc O_{\lambda})$, then there is a natural isomorphism
\[
\Lie(G) \otimes_{\ZZ} \FF_{\ell} \to G(\mc O/\lambda^n)_{n-1}, \qquad A \mapsto I_d + \lambda^{n-1} \cdot A.
\]
We want now to prove a similar isomorphism for $\SU(V/\lambda^n, \bb e)_{n-1}$. It turns
out however that since $\SU$ is not an algebraic group over $\mc O_{\lambda}$, the result will actually depend on the parity of~$n$.
Denote for any $n \ge 1$ and any ring~$R$:
\[
\mf{su}^{(n)}(R) := \{ A \in M_d(R) : \Gamma A = (-1)^n \cdot A^T \Gamma, \, \, \tr(A) = 0 \}.
\]
The following lemma shows that the spaces $\mf{su}^{(n)}(\FF_{\ell})$ play the role of the Lie algebra $\Lie(G)$ for the unitary group.
\begin{lemma} \label{lem:lie_of_su}
	The map:
	\begin{align*}
		\mf{su}^{(n)}(\FF_{\ell}) &\to \Gl_d(\mc O/\lambda^n)_{n-1}\\
		A &\mapsto I_d + \lambda^{n-1} A
	\end{align*}
	is an isomorphism onto $\SU(V/\lambda^n, \bb e)_{n-1}$.
\end{lemma}
\begin{proof}
	Suppose that $\bb A := I_d + \lambda^{n-1} A_{n-1} \in \Gl_d(\mc O/\lambda^n)_{n-1}$.
	The congruence~\eqref{eqn:ol_lambda_is_minus_lambda} implies that in $\Gl_d(\mc O/\lambda^n)$:
	\begin{align*}
		\bb A^{\dagger} \Gamma \bb A &= (I_d + \ol{\lambda}^{n-1} A_{n-1}^T) \cdot \Gamma \cdot (I_d + \lambda^{n-1} A_{n-1})\\
		&= \Gamma + \lambda^{n-1} \cdot ((-1)^{n-1} A_{n-1}^T \Gamma + \Gamma A_{n-1}).
	\end{align*}
	This shows that $\bb A^{\dagger} \Gamma \bb A = \Gamma$ if and only if $\Gamma A_{n-1} = (-1)^n A_{n-1}^T \Gamma$. Moreover,
	the equality~\eqref{eqn:trace_and_det} implies that $\det_{\mc O_{\lambda}} \bb A = 1$ if and only if $\tr(A_{n-1}) = 0$.
\end{proof}
\noindent For any $1 \le i, j \le d$ and $n \ge 1$ define the following matrices:
\[
E_{ij}^{(n)} := 
E_{ij} + (-1)^n \cdot E_{ji},
\]
where $E_{ij}$ is the matrix with one on the position $(i, j)$ and zeros on other positions.
In the sequel, we use the following equalities for any $1 \le i, j, l \le d$, $j \neq i, l$:
\begin{align}
	[\Gamma^{-1} E_{ij}^{(m)}, \Gamma^{-1} E_{jl}^{(n)}] &= 
	\begin{cases}
		\alpha_j^{-1} \cdot \Gamma^{-1} E^{(m+n+1)}_{il}, & i \neq l, \\
		(1 + (-1)^{n+m+1}) \cdot \alpha_i^{-1} \cdot \alpha_j^{-1} \cdot (E_{ii} - E_{jj}), & i = l.
	\end{cases}
	\label{eqn:Eij_commutators}
\end{align}
In fact, they follow from the equalities $[E_{ij}, E_{kl}] = \delta_{jk} E_{il} - \delta_{li} E_{kj}$
and $\Gamma^{-1} E_{ij} = \alpha_i^{-1} E_{ij}$.
\newpage
\begin{lemma} \label{lem:basis_of_su_n}
	\begin{enumerate}[(1)]
		\item[]
		\item If $2 \nmid n$, then the matrices $\Gamma^{-1} E_{ij}^{(n)}$, $1 \le i < j \le d$ are the basis of the $\FF_{\ell}$-vector space $\mf{su}^{(n)}(\FF_{\ell})$. In particular
		$\dim_{\FF_{\ell}} \mf{su}^{(n)}(\FF_{\ell}) = {d \choose 2}$.
		
		\item If $2 | n$, then the matrices $\Gamma^{-1} E_{ij}^{(n)}$, $1 \le i < j \le d$ and the matrices $E_{ii} - E_{dd}$, $1 \le i < d$, are the basis of the $\FF_{\ell}$-vector space $\mf{su}^{(n)}(\FF_{\ell})$. In particular
		$\dim_{\FF_{\ell}} \mf{su}^{(n)}(\FF_{\ell}) = {d + 1 \choose 2} - 1$.
	\end{enumerate}
\end{lemma}
\begin{proof}
	\begin{enumerate}[(1),leftmargin=*]
		\item[]
		
		\item This follows easily from the fact that the map $A \mapsto \Gamma^{-1} A$ is an isomorphism between the space of antisymmetric matrices in $M_d(\FF_{\ell})$
		and $\mf{su}^{(n)}(\FF_{\ell})$.
		
		\item As in (1), note that we have an isomorphism $A \mapsto \Gamma^{-1} A$ between the space of symmetric matrices in $M_d(\FF_{\ell})$
		and $\{ A : \Gamma A = A^T \Gamma \}$. Thus the basis of the latter space is given by the matrices $\Gamma^{-1} E_{ij}^{(n)}$ for $1 \le i \le j \le d$.
		By adding the condition on the vanishing of the trace we easily see that the basis of $\mf{su}^{(n)}(\FF_{\ell})$ is
		\[
		\Gamma^{-1} E_{ij} - \frac{\tr(\Gamma^{-1} E_{ij})}{\tr(\Gamma^{-1} E_{dd})} \cdot \Gamma^{-1} E_{dd}
		=
		\begin{cases}
			\Gamma^{-1} E_{ij}, & i \neq j\\
			\alpha_i^{-1} \cdot (E_{ii} - E_{dd}), & i = j,
		\end{cases}
		\]
		for $1 \le i \le j < d$, $(i, j) \neq (d, d)$.
	\end{enumerate}
\end{proof}
\begin{corollary} \label{cor:X=GY+YG}
	For every $X \in M_d(\FF_{\ell})$ such that $X^T = (-1)^n \cdot X$ there exists $Y \in M_d(\FF_{\ell})$
	such that $X = \Gamma Y + (-1)^{n-1} \cdot Y^T \Gamma$. Moreover:
	\begin{itemize}
		\item if $2 | n$, then the matrix $Y$ can be chosen to have an arbitrary trace $c \in \FF_{\ell}$,
		
		\item if $2 \nmid n$, then $\tr(Y) = \frac 12 \tr(\Gamma^{-1} X)$.
	\end{itemize}
\end{corollary}
\begin{proof}
	Fix $n$ and consider the map
	\begin{align*}
		\varphi : M_d(\FF_{\ell}) &\to M_d(\FF_{\ell}) \times \FF_{\ell},\\
		Y &\mapsto (\Gamma Y + (-1)^{n-1} Y^T \Gamma, \tr(Y)).
	\end{align*}
	Note that its kernel is $\mf{su}^{(n)}(\FF_{\ell})$ and that $\im \varphi \subset \mc Z$, where
	\begin{equation}
		\mc Z :=
		\begin{cases}
			\{ A : A^T = -A \} \times \FF_{\ell}, & \textrm{ if } 2|n,\\
			\{ (A, \frac{1}{2} \tr(\Gamma^{-1} A)) : A^T = A \}, & \textrm{ if } 2\nmid n.
		\end{cases}
		\label{eqn:im_phi}
	\end{equation}
	One easily computes that
	\[
	\dim_{\FF_{\ell}} \mc Z =
	\begin{cases}
		{d \choose 2} + 1, & 2 | n,\\
		{d+1 \choose 2}, & 2 \nmid n.	
	\end{cases}
	\]
	On the other hand, one checks that $\dim_{\FF_{\ell}} \textrm{im} \varphi = d^2 - \dim_{\FF_{\ell}} \mf{su}^{(n)}(\FF_{\ell}) = \dim_{\FF_{\ell}} \mc Z$
	by using Lemma~\ref{lem:basis_of_su_n}. Thus $\im \varphi = \mc Z$, which ends the proof.
\end{proof}
\begin{lemma} \label{lem:su_surjection}
	Let $n \ge 2$. Then the natural map $\SU(V/\lambda^n, \bb e)_1 \to \SU(V/\lambda^{n - 1}, \bb e)_1$
	is a surjection.
\end{lemma}
\begin{proof}
	Let $\bb A \in \SU(V/\lambda^{n - 1}, \bb e)_1$. Take any lift $\bb A' \in \Gl_d(\mc O/\lambda^n)$ of $\bb A$
	and suppose that $(\bb A')^{\dagger} \Gamma \bb A' = \Gamma + \lambda^{n - 1} \cdot X$ for some $X \in M_d(\FF_{\ell})$
	and $\det_{\mc O_{\lambda}} \bb A' = 1 + \lambda^{n-1} \cdot c$ for some $c \in \FF_{\ell}$. By~\eqref{eqn:trace_and_det} it is easy to check that
	$\tr \Gamma^{-1} X = c + (-1)^{n-1} \cdot c$.
	Note that $(\bb A')^{\dagger} \Gamma \bb A'$ is Hermitian, so that by~\eqref{eqn:ol_lambda_is_minus_lambda} we obtain:
	\[
	(\Gamma + \lambda^{n - 1} \cdot X)^{\dagger} = \Gamma + \lambda^{n - 1} \cdot X \Rightarrow X^T = (-1)^{n-1} \cdot X.
	\]
	By Corollary~\ref{cor:X=GY+YG} $X$ is of the form $\Gamma Y + (-1)^{n-1} \cdot Y^T \Gamma$ for some $Y \in M_d(\FF_{\ell})$.
	Moreover, if $2 | n$, we can assume that $\tr(Y) = c$. If $2 \nmid n$, this is automatic, since then $c = \frac{1}{2} \tr(\Gamma^{-1} X) = \tr(Y)$. Let $\bb A'' := \bb A' - \lambda^{n-1} Y$. Then:
	\begin{align*}
		(\bb A'')^{\dagger} \cdot \bb A'' = \Gamma + \lambda^{n - 1} \cdot X - \lambda^{n-1} \cdot (Y \Gamma + (-1)^{n-1} \cdot Y^T \Gamma) = \Gamma.
	\end{align*}
	Moreover, using the fact that $\bb A' \equiv I_d \pmod{\lambda}$ and~\eqref{eqn:trace_and_det}, we have in $\mc O/\lambda^n$:
	\begin{align*}
		\det {}_{\mc O_{\lambda}} \bb A'' &= \det {}_{\mc O_{\lambda}}(\bb A' - \lambda^{n-1} \cdot Y) = \det {}_{\mc O_{\lambda}}(\bb A') \cdot \det {}_{\mc O_{\lambda}}(I_d - \lambda^{n-1} \cdot (\bb A')^{-1} Y)\\
		&= \det {}_{\mc O_{\lambda}}(\bb A') \cdot (1 - \lambda^{n-1} \tr((\bb A')^{-1} Y))
		= \det {}_{\mc O_{\lambda}}(\bb A') \cdot (1 - \lambda^{n-1} \cdot \tr(Y))\\
		&= (1 + \lambda^{n-1} \cdot c) (1 - \lambda^{n-1} \cdot c) = 1.
	\end{align*}
	Hence $\bb A''$ is a lift of $\bb A$ that belongs to $\SU(V/\lambda^n, \bb e)$.
\end{proof}
\begin{corollary} \label{cor:su_inside_commutator}
	Let $n \ge 3$, $d \ge 3$, $N := \lfloor \frac{n - 1}{2} \rfloor$, $M := n - 1 - N$. Then:
	\[
 \llbracket \SU(V/\lambda^n, \bb e)_N, \SU(V/\lambda^n, \bb e)_M \rrbracket = 
		\SU(V/\lambda^n, \bb e)_{n-1}.
	\]
\end{corollary}
\begin{proof}
	Note that $M = N$, if $n$ is odd and $M = N+1$ if $n$ is even.
	Therefore by Lemma~\ref{lem:computation_of_commutator},
	for any $\bb A \in \SU(V/\lambda^n, \bb e)_N$, $\bb B \in \SU(V/\lambda^n, \bb e)_M$ we have 
	\[
  \llbracket \bb A, \bb B \rrbracket = I + \lambda^{n-1} \cdot [A_N, B_M].
	\]
	This equality implies the inclusion:
	\begin{align*}
		\llbracket \SU(V/\lambda^n, \bb e)_N, \SU(V/\lambda^n, \bb e)_M \rrbracket &\subset 
		\SU(V/\lambda^n, \bb e)_{n-1}.
	\end{align*}
	By Lemmas~\ref{lem:lie_of_su} and~\ref{lem:basis_of_su_n} we have $I_d + \lambda^N \cdot \Gamma^{-1} E_{ij}^{(N+1)} \in \SU(V/\lambda^{N+1}, \bb e)_N$, $I_d + \lambda^M \cdot \Gamma^{-1} E_{jl}^{(M+1)} \in \SU(V/\lambda^{M+1}, \bb e)_M$. 
	Let $\bb A$ be any lift of $I_d + \lambda^N \cdot \Gamma^{-1} E_{ij}^{(N+1)}$ to $\SU(V/\lambda^n, \bb e)$ 
	and let $\bb B$ be any lift of $I_d + \lambda^M \cdot \Gamma^{-1} E_{jl}^{(M+1)}$ to $\SU(V/\lambda^n, \bb e)$ (we can pick $\bb A$ and $\bb B$ using Lemma~\ref{lem:su_surjection}). Then, by Lemma~\ref{lem:computation_of_commutator}
	and by~\eqref{eqn:Eij_commutators}:
	\begin{align*}
		\llbracket \bb A, \bb B \rrbracket &= I_d + \lambda^{n-1} \cdot [\Gamma^{-1} E_{ij}^{(N+1)}, \Gamma^{-1} E_{jl}^{(M+1)}]\\
		&=
		\begin{cases}
			I_d + \lambda^{n-1} \cdot \alpha_j^{-1} \cdot \Gamma^{-1} E_{il}^{(n)}, & i \neq l,\\
			I_d + \lambda^{n-1} \cdot (1 + (-1)^n) \cdot \alpha_i^{-1} \cdot \alpha_j^{-1} \cdot (E_{ii} - E_{jj}), & i = l.
		\end{cases}
	\end{align*}
	Now it suffices to observe that the matrices $\Gamma^{-1} E_{il}^{(n)}$ and, if $2|n$, $E_{ii} - E_{dd}$ generate $\mf{su}^{(n)}(\FF_{\ell})$ by Lemma~\ref{lem:basis_of_su_n} and to use
	Lemma~\ref{lem:lie_of_su}.
\end{proof}

\begin{proof}[Proof of Proposition~\ref{prop:lifting}]
	Since $\mc H$ is closed, it suffices to show that
	$\mc H$ surjects onto $\SU(V/\lambda^n, \bb e)^{\mc S}$ for every $n \ge 2$.
	We prove this by induction on $n$. For $n = 2$ this follows by the assumption.
	Suppose now that $\mc H$ surjects onto $\SU(V/\lambda^{n-1}, \bb e)^{\mc S}$.
	By abuse of notation, write $\mc H$ for the image of $\mc H$ in $\SU(V/\lambda^n, \bb e)$.
	Let also for $i < n$:
	\[
	\mc H_i := \ker(\mc H \to \SU(V/\lambda^i, \bb e)).
	\]
	It suffices to show that
	\[
	\SU(V/\lambda^n, \bb e)_{n-1} \le \mc H.
	\]
	Let $N$ and $M$ be as in Corollary~\ref{cor:su_inside_commutator}.
	Note that $\mc H_{n-1}$ is a normal subgroup of $\SU(V/\lambda^n, \bb e)_1$, since it is a subgroup of the central subgroup $\SU(V/\lambda^n, \bb e)_{n-1}$
	(cf. Lemma~\ref{lem:gl_in_center}). Moreover, since by assumption for any $i \ge 1$ the group $\mc H_i$ surjects onto
	$\SU(V/\lambda^{n - 1}, \bb e)_i$ we have:
	\[
	\mc H_i/\mc H_{n-1} \cong \SU(V/\lambda^{n - 1}, \bb e)_i.
	\]
	Thus the exact sequence:
	\[
	0 \to \SU(V/\lambda^n, \bb e)_{n-1}/\mc H_{n-1} \to \SU(V/\lambda^n, \bb e)_i/\mc H_{n-1} \to \SU(V/\lambda^{n - 1}, \bb e)_i \to 0
	\]
	has a section, mapping $\SU(V/\lambda^{n - 1}, \bb e)_i$ onto $\mc H_i/\mc H_{n-1} \le \SU(V/\lambda^n, \bb e)_i/\mc H_{n-1}$.
	Hence:
	\begin{equation} \label{eqn:SU/H_is_a_product}
		\SU(V/\lambda^n, \bb e)_i/\mc H_{n-1} \cong \SU(V/\lambda^{n - 1}, \bb e)_i \times \bigg( \SU(V/\lambda^n, \bb e)_{n-1}/\mc H_{n-1} \bigg)
	\end{equation}
	(it is a direct product, since $\SU(V/\lambda^n, \bb e)_{n-1}$ is in the center of $\SU(V/\lambda^{n - 1}, \bb e)_i$).
	Observe that $\SU(V/\lambda^{n-1}, \bb e)_M$ is
	in the center of $\SU(V/\lambda^{n-1}, \bb e)_N$
	by Lemma~\ref{lem:gl_in_center}, since $M + N = n-1$.
	Therefore (using~\eqref{eqn:SU/H_is_a_product} for $i = M, N$) $\SU(V/\lambda^n, \bb e)_M/\mc H_{n-1}$ is in the center of $\SU(V/\lambda^n, \bb e)_N/\mc H_{n-1}$. This yields:
	\[
		\llbracket \SU(V/\lambda^n, \bb e)_N, \SU(V/\lambda^n, \bb e)_M \rrbracket \subset \mc H_{n-1}.
	\]
	Hence, by Corollary~\ref{cor:su_inside_commutator} we have $\SU(V/\lambda^n, \bb e)_{n-1} \subset \mc H_{n-1}$, which ends the proof.
\end{proof}

\section{Proof of Theorem~\ref{thm:image_of_rho}} \label{sec:proof_of_thm_rho}
In order to apply Proposition~\ref{prop:lifting}, we need additional information on the $\lambda^2$-torsion field
of $J$. By descent theory, see e.g. \cite[Corollary~3.0.8.]{arul_torsion_fermat_quotients}:
\begin{equation*}
	K(J[\lambda^2]) =
	K(\alpha_1, \ldots, \alpha_r, \sqrt[\ell]{\alpha_i - \alpha_j} : 1 \le i < j \le r).
\end{equation*}
Thus $[K(J[\lambda^2]) : K(J[\lambda])] \le \ell^{{r \choose 2}}$. The following lemma shows that under some mild assumptions
the equality is achieved.
\begin{lemma} \label{lem:l_independence}
	Let $f \in L[x]$ be a monic separable polynomial of degree $r$ over a number field~$L$. Let $\{ \alpha_1, \ldots, \alpha_r \}$ be the roots of $f$ and let $M$ be its splitting field. Suppose that:
	\begin{enumerate}[(1)]
		\item $\Gal(f)$ acts $2$-transitively on $\{ \alpha_1, \ldots, \alpha_r  \}$,
		
		\item that there exists
		a prime $\mf p \nmid 2 \disc(L/\QQ)$ such that $\ord_{\mf p}(\disc(f)) = 1$ and that the coefficients of $f$ are $\mf p$-integral.
	\end{enumerate}
	Then the elements $\{ \alpha_i - \alpha_j : 1 \le i < j \le r \}$ are independent in the $\FF_{\ell}$-vector space $L^{\times}/(L^{\times})^{\ell}$.
\end{lemma}
\begin{proof}
	Denote by $\mc P$ any prime ideal of $\mc O_M$ lying over~$\mf p$.
	Since $\ord_{\mf p}(\disc(f)) = 1$ and $\mf p \nmid 2 \disc(L/\QQ)$,
	the polynomial $f$ has a single double root modulo $\mf p$
	(see e.g. \cite[Lemma~2.1]{Garnek_maximally_disjoint} for the case $L = \QQ$; the general case can be proven in a similar manner). Thus $\mc P | \alpha_i - \alpha_j$ for precisely
	one pair $1 \le i < j \le r$, without loss of generality $\mc P | \alpha_1 - \alpha_2$. 
	We claim that $e(\mc P/\mf p)$ (the ramification index of $\mc P$ over $\mf p$) equals~$2$.
	Indeed, let $\mc P_1 := \mc P \cap \mc O_{L(\alpha_1)}$ and let $f = (x - \alpha_1) \cdot f_1$ for some $f_1 \in L(\alpha_1)[x]$. Then $\mc P_1 \nmid \disc(f_1) = \prod_{1 < i < j \le r} (\alpha_i - \alpha_j)^2$. Therefore $\mc P_1$ is unramified in $M$, the splitting field of $f_1$, and $e(\mc P/\mf p) = e(\mc P_1/\mf p)$. Moreover, by \cite[Theorem~3.41 and Remark~3.43]{milneANT} we have $e(\mc P_1/\mf p) = 2$, which proves the claim. Thus:
	\[
	\ord_{\mc P}(\alpha_1 - \alpha_2) \le \frac 12 \ord_{\mc P}(\disc(f)) = \ord_{\mf p}(\disc(f)) = 1,
	\]
	i.e. $\ord_{\mc P}(\alpha_1 - \alpha_2) = 1$. 
	Since $\Gal(f)$ is a $2$-transitive subgroup of $S_r$, 
	for any $1 \le i < j \le r$ we may choose $\sigma_{ij} \in \Gal(f)$ such that
	$\sigma_{ij}(\alpha_1) = \alpha_i$, $\sigma_{ij}(\alpha_2) = \alpha_j$.
	Define $\mc P_{ij} := \sigma_{ij}(\mc P)$. Then:
	\[
	\ord_{\mc P_{ij}}(\alpha_k - \alpha_l) =
	\begin{cases}
		1, & (k, l) = (i, j),\\
		0, & (k, l) \neq (i, j).
	\end{cases}
	\]
	Suppose now that for some $a_{ij} \in \ZZ$ we have:
	\begin{equation*}
		\prod_{i < j} (\alpha_i - \alpha_j)^{a_{ij}} = \gamma^{\ell} \qquad \textrm{ for some } \gamma \in M^{\times}.
	\end{equation*}
	Then, by computing the valuation of both sides with respect to $\mc P_{ij}$, we obtain $a_{ij} = \ell \cdot \ord_{\mc P_{ij}}(\gamma)$,
	i.e. $\ell | a_{ij}$. This ends the proof. 
\end{proof}
\begin{proof}[Proof of Theorem~\ref{thm:image_of_rho}]
	From the discussion in Section~\ref{sec:preliminaries} it follows that $\Gal(K(J[\lambda^2])/K)$ is contained in $\U_{r-1}^{\varepsilon}(\mc O/\lambda^2)^{\Gal(f)}$. On the other hand, one easily proves (as in Lemmas~\ref{lem:lie_of_su} and \ref{lem:basis_of_su_n}) that $\U_{r-1}^{\varepsilon}(\mc O/\lambda^2)_1$
	is isomorphic to the space of symmetric $(r-1) \times (r-1)$-matrices and thus $|\U_{r-1}^{\varepsilon}(\mc O/\lambda^2)_1| = \ell^{{r \choose 2}}$. Hence by Lemma~\ref{lem:l_independence} and by Kummer theory we have:
	\begin{align*}
		|\Gal(K(J[\lambda^2])/K)| &= [K(J[\lambda^2]):K(J[\lambda])] \cdot [K(J[\lambda]):K]
		= \ell^{{r \choose 2}} \cdot |\Gal(f)|\\
		&= |\U_{r-1}^{\varepsilon}(\mc O/\lambda^2)_1| \cdot |\Gal(f)| =  |\U_{r-1}^{\varepsilon}(\mc O/\lambda^2)^{\Gal(f)}|,
	\end{align*}
	which implies that $\Gal(K(J[\lambda^2])/K) = \U_{r-1}^{\varepsilon}(\mc O/\lambda^2)^{\Gal(f)}$.
	
	Let $\mc H := \ker(\det_{\mc O_{\lambda}} : \rho_{J, \ell}(\Gal_K) \to \mc O_{\lambda}^{\times})$. Note that $\mu(\mc H) \subset \chi_{cyc}(\Gal_K) = 1 + \ell \ZZ_{\ell}$. On the other hand, $\mu(\mc H) \subset \mu_{r-1}$ by~\eqref{eqn:ol_det_det_equals_mu_r-1}.
	Hence $\mu(\mc H) = 1$ and $\mc H \subset \SU^{\varepsilon}_{r-1}(\mc O_{\lambda})$.
	Moreover, $S_r \cap \SU_{r-1}(\FF_{\ell}) = A_r$.
	Therefore the reduction of $\mc H$ modulo $\lambda^2$ equals 
	\[
		\Gal(K(J[\lambda^2])/K) \cap \SU_{r-1}^{\varepsilon}(\mc O/\lambda^2)^{\Gal(f) \cap A_r} = \SU_{r-1}^{\varepsilon}(\mc O/\lambda^2)^{\Gal(f)\cap A_r}.
	\]
	Thus by Proposition~\ref{prop:lifting}
	we obtain $\mc H = \SU_{r-1}^{\varepsilon}(\mc O_{\lambda})^{\Gal(f) \cap A_r}$. Since the map $\det_{\mc O_{\lambda}} : \im \rho_{J, \ell} \to \mc D_J$
	is surjective (by definition), one sees that $\im \rho_{J, \ell} = \GU_{r-1}^{\varepsilon}(\mc O_{\lambda})^{\Gal(f)}_{\det \in \mc D_J}$.
\end{proof}

\section{The endomorphism character} \label{sec:endo_char}
In the previous section we determined the image of the $\lambda$-adic representation of~$J$
up to the image of the determinant (denoted $\mc D_J$). In this section
we prove some results on $\mc D_J$. Note that $\mc D_J \subset \mc U_{\ell, r}$. Indeed, this follows from~\eqref{eqn:ol_det_det_is_cyc} and from the equality $\chi_{cyc}^{r-1}(\Gal_K) = (1 + \ell \ZZ_{\ell})^{r-1} = 1 + \ell \cdot (r-1) \ZZ_{\ell}$. We prove now a result concerning the torsion in $\mc D_J$, which is an easy consequence of previous observations. 
\begin{proposition} \label{prop:torsion_in_DJ}
	If $J$ satisfies the assumptions~(1) and~(2) of Theorem~\ref{thm:image_of_rho}, then $\mu_{2 \ell} \subset \mc D_J$.
\end{proposition}
\begin{proof}
	From the proof of Theorem~\ref{thm:image_of_rho} in the previous section it follows that
	$\Gal(K(J[\lambda^2])/K) = \U^{\varepsilon}_{r-1}(\mc O/\lambda^2)^{\Gal(f)}$. In particular,
	$\mc D_J$ surjects onto the reduction of $\mc U_{\ell, r}$ modulo $\lambda^2$.
	But this reduction equals $\mu_{2 \ell}$, which ends the proof.
\end{proof}
\begin{remark} \label{rem:DJ_and_mu_2}
	Note that the composition $S_r \to \Gl_{r-1}(\FF_{\ell}) \stackrel{\det}{\longrightarrow} \FF_{\ell}^{\times}$ is the signature homomorphism. Thus the reduction of $\mc D_J$
	modulo $\lambda$ contains $\mu_2$ if and only if $\Gal(f) \le S_r$ is not a subgroup of $A_r$
	(the alternating group of degree $r$).
\end{remark}
To say more about~$\mc D_J$, we need the theory of Hecke characters. Let $L$ be an arbitrary number field. Consider the the group of ideles of $L$, i.e.
\[
\II_L := \{ (a_v)_v \in \prod_v L_v^{\times} : a_v \in \mc O_{L, v}^{\times} \textrm{ for almost all } v \}
\]
where the product is taken over all (finite and infinite) places of $L$. Recall that the class field
theory yields an epimorphism $\II_L/L^{\times} \to \Gal_L^{ab}$. A continuous homomorphism
$\chi : \II_L/L^{\times} \to \CC^{\times}$ is called a \bb{Hecke character}.
For any place $v$, denote by $\iota_v : L_v^{\times} \hookrightarrow \II_L$ the natural embedding.
For every finite $v$ there exists a number $m_v \ge 0$ such that $\chi(\iota_v(1 + \mf p_v^{m_v})) = \{ 1 \}$
(where we write $1 + \mf p_v^0 := \mc O_v^{\times}$). If $m_v$ can be taken to be $0$, we say that $\chi$ is \bb{unramified} at~$v$. In particular, for any $a \in \mc O_{L, v}^{\times}$, $\chi(\iota_v(a))$ is a root of unity, since $a^n \in 1 + \mf p_v^{m_v}$
for $n$ large enough. Let $L_{\infty}^{\times} := \prod_{v | \infty} L_v^{\times}$ and write $L_{\infty}^{\times, \circ}$ for the connected component
of $L_{\infty}^{\times}$ (note that if $L_{\infty}^{\times} = (\CC^{\times})^{r_1} \times (\RR^{\times})^{r_2}$, then $L_{\infty}^{\times, \circ} = (\CC^{\times})^{r_1} \times (\RR_+)^{r_2}$).
We say that a character $L_{\infty}^{\times, \circ} \to \CC^{\times}$ is \bb{algebraic}, if it comes from a homomorphism $T : L^{\times} \to \CC^{\times}$,
$T = \prod_{\sigma : L \hookrightarrow \CC} \sigma^{n_{\sigma}}$ for some integers $n_{\sigma}$. We say that a Hecke character
$\chi$ is \bb{algebraic}, if $\chi|_{L_{\infty}^{\times, \circ}}$ is algebraic. In this case, we say that the corresponding homomorphism $T$ is the \bb{infinity type} of~$\chi$. If $\chi$ is
algebraic, its image must lie in a number field $M$. Fix a rational prime $p$ and
let~$\mf q$ be a prime of $M$ over $p$. The homomorphism $T$ extends uniquely
to a homomorphism $T_{\mf q} : \prod_{v | p} L_v^{\times} \to M_{\mf q}^{\times}$.
This allows one to define a character $\chi_{\mf q} : \II_L \to M_{\mf q}^{\times}$
(the \bb{$\mf q$-adic avatar of $\chi$}) by the formula:
\begin{equation} \label{eqn:avatar}
	\chi_{\mf q}(a) := \chi(a) \cdot T(a_{\infty})^{-1} \cdot T_{\mf q}(a_p),
\end{equation}
where $a_p$ denotes the components of $a$ lying in $\prod_{v | p} L_v^{\times}$. One shows that $\chi_{\mf q}$ descends to a Galois character $\Gal_L \to M_{\mf q}^{\times}$, which we also denote by $\chi_{\mf q}$.\\ 

We return now to Setup~\ref{setup}. Consider the
\bb{endomorphism character of $J$}, i.e. the Galois character $\Omega_{J, \lambda} : \Gal_K \to K_{\lambda}^{\times}$, $\Omega_{J, \lambda} := \det_{\mc O_{\lambda}} \circ \rho_{J, \ell}$.
The endomorphism character has been studied mostly in case when the considered abelian variety has complex multiplication
-- for example the article \cite{Coleman_stable_reductions} relates the endomorphism character of the Jacobian of the Fermat curve to
Jacobi sums. The article \cite{goodman_superelliptic} studies the reduction of the endomorphism character of a superelliptic Jacobian
modulo a prime different then~$\lambda$.
As explained in \cite{goodman_superelliptic}, the endomorphism character comes from an algebraic Hecke character $\Omega_J : \II_K/K^{\times} \to K^{\times}$. 
Moreover, the infinity type of $\Omega_J$ can be computed knowing the action of $\ZZ/\ell$ on $H^0(C, \Omega_{C/\CC})$, see \cite[Proposition~14]{fite_ordinary_primes}. An explicit formula for this action in case of superelliptic
curves in given e.g. in \cite[Remark~3.7]{Zarhin_endo_Jacobians_cyclic_covers}. Thus the 
infinity type of $\Omega_J$ is given by the formula:
\[
T_{\ell, r} := \prod_{j = 1}^{\ell - 1} \sigma_j^{n(j)},
\]
where $\sigma_j \in \Gal(K/\QQ)$ is defined by $\sigma_j(\zeta_{\ell}) = \zeta_{\ell}^j$ and $n(j) := \lfloor \frac{r \cdot (\ell - j)}{\ell} \rfloor$ (cf. \cite[\S 5.3]{goodman_superelliptic}). Define $\mc T_{\ell, r} := T_{\ell, r}(\mc O_{\lambda}^{\times}) \subset \mc O_{\lambda}^{\times}$.
We show below (Proposition~\ref{prop:Tlr_in_DJ}) that $[\mc U_{\ell, r} : \mc D_J]$ is finite if and only if $[\mc U_{\ell, r} : \mc T_{\ell, r}]$ is finite.
\begin{example}
	Let $C$ be the curve $y^3 = x^2 - D$, where $D \in \ZZ$, $D \neq 0$. 
	Then $\mc O = \ZZ[\zeta_3]$. By \cite[Example II.10.6]{Silverman_Advanced},
	$\Omega_J$ is unramified for primes $\mf p \nmid 6 D$ and for any such prime $\mf p$
	with uniformizer $\pi_{\mf p}$, $\pi_{\mf p} \equiv 2 \pmod 3$, we have:
	\[
	\Omega_J(\pi_{\mf p}) = - \ol{\prs{4D}{\pi_{\mf p}}{6}} \cdot \pi_{\mf p}.
	\]
	Note that these conditions determine $\Omega_J$ uniquely. Moreover,
	the above formula shows that the infinity type of $\Omega_J$ is $T = \id$, which is consistent with the above formula.
\end{example}
\begin{proposition} \label{prop:Tlr_in_DJ}
	We have $\mc T_{\ell, r} \subset \mu_{2 \ell} \mc D_J$. Moreover:
	\[
	[\mu_{2 \ell} \mc D_J : \mu_{2 \ell} \mc T_{\ell, r}] \le \ell^{\ord_{\ell}(h_{\ell})},
	\]
	where $h_{\ell}$ is the class number of $\QQ(\zeta_{\ell})$.
\end{proposition}
\begin{proof}
	We first prove that $\mc T_{\ell, r} \subset \mu_{2 \ell} \mc D_J$. 
	Indeed, let $a \in \mc O_{\lambda}^{\times}$. Recall that $\Omega_J(\iota_{\lambda}(a))$ is a root of unity
	in $K$. Therefore $\Omega_J(\iota_{\lambda}(a)) \in \mu_{2 \ell}$ and by~\eqref{eqn:avatar}:
	\begin{align*}
		T_{\ell, r}(a) =  \Omega_J(\iota_{\lambda}(a))^{-1} \cdot \Omega_{J, \lambda}(\iota_{\lambda}(a)) \in \mu_{2 \ell} \mc D_J.
	\end{align*}
	This shows the desired inclusion. We prove now the bound on the index of $\mu_{2 \ell} \mc T_{\ell, r}$ in $\mu_{2 \ell} \mc D_J$. Observe that the group $\mc O_{\lambda}^{\times}/\mu_{\ell - 1}$ is a pro-$\ell$ group and that $\mu_{\ell - 1} \cap \mc U_{\ell, r} = \mu_2$. Therefore $\mc U_{\ell, r}/\mu_2$
	is also a pro-$\ell$ group and it is enough to show that $[\mu_{2 \ell} \mc D_J : \mu_{2 \ell} \mc T_{\ell, r}]$ divides $h_{\ell}$. Clearly, $\Omega_{J, \lambda}(K_v^{\times}) = \{ 1 \}$ for any infinite place $v$. For any finite place $v$ and $a \in \mc O_v^{\times}$ we have either
	$\Omega_{J, \lambda}(\iota_v(a)) = \Omega_J(\iota_v(a)) \in \mu_{2 \ell}$ if $v \neq \lambda$ or $\Omega_{J, \lambda}(\iota_v(a)) = \Omega_J(\iota_v(a)) \cdot T_{\ell, r}(a) \in \mu_{2 \ell} \mc T_{\ell, r}$ if $v = \lambda$.
	Moreover, since $\lambda \in \mc O_v^{\times}$ for $v \neq \lambda$,
	we have
	\[
	\Omega_{J, \lambda}(\iota_{\lambda}(\lambda)) = \prod_{v \neq \lambda}
	\Omega_{J, \lambda}(\iota_v(\lambda))^{-1} \in \mu_{2 \ell}.
	\]
	Finally, if $v \neq \lambda$ is a finite place, then $\mf p_v^{h_{\ell}} = \alpha \mc O$ is a principal ideal. Thus $\pi_v^{h_{\ell}} = u_1 \cdot \alpha$ for some $u_1 \in \mc O^{\times}_v$. Observe also that $\alpha \in \mc O_{v'}^{\times}$ for $v' \neq v$.
	Hence:
	\begin{align*}
		\Omega_{J, \lambda}(\pi_v)^{h_{\ell}} &= 
		\Omega_{J, \lambda}(\iota_v(u_1)) \cdot \Omega_{J, \lambda}(\iota_v(\alpha))\\
		&= \Omega_{J, \lambda}(\iota_v(u_1)) \cdot \prod_{v' \neq v} \Omega_{J, \lambda}(\iota_{v'}(\alpha))^{-1}\\
		&= \left(\Omega_{J, \lambda}(\iota_v(u_1)) \cdot \prod_{v' \neq v, \lambda} \Omega_{J, \lambda}(\iota_{v'}(\alpha))^{-1} \cdot \Omega_J(\iota_{\lambda}(\alpha))^{-1} \right) \cdot T_{\ell, r}(\alpha)^{-1}  \in \mu_{2 \ell} \mc T_{\ell, r}.
	\end{align*}
	This ends the proof.
\end{proof}
\noindent Note that $\mc O_{\lambda}^{\times} \cong \mu_{\ell (\ell - 1)} \times (1 + \lambda^2 \mc O_{\lambda})$ and
\[
\mc U_{\ell, r} = \mu_{2 \ell} \times (1 + \ell \cdot (r-1) \ZZ_{\ell}) \times \mc U_{\ell}',
\]
where $\mc U_{\ell}' := \{ d \in 1 + \lambda^2 \mc O_{\lambda} : d \cdot \ol d = 1 \}$.
Moreover, $T_{\ell, r} = \prod_{j = 1}^{\ell - 1} \sigma_j^{(r-1)/2} \cdot T_{\ell, r}'$, where
$T_{\ell, r}' := \prod_{j = 1}^{\ell - 1}  \sigma_j^{n'(j)}$ and
\begin{equation} \label{eqn:n_prim_eqn}
	n' : (\ZZ/\ell)^{\times} \to \frac 12 \ZZ, \quad n'(j) :=  \left\lfloor \frac{(\ell - j) r}{\ell} \right\rfloor - \frac{r-1}{2}
\end{equation}
(here we identify $(\ZZ/\ell)^{\times}$ with $\{ 1, \ldots, \ell - 1 \}$).
Observe that
\begin{equation} \label{eqn:n'_is_odd}
	n'(j) := \frac{r-1}{2} - \left\lfloor \frac{j r}{\ell} \right\rfloor
	= -n'(\ell-j).
\end{equation}
This implies in particular that $\mc T'_{\ell, r} := T_{\ell, r}'(1 + \lambda^2 \mc O_{\lambda})$ is contained
in $\mc U_{\ell}'$. Moreover, $(\prod_{j = 1}^{\ell - 1} \sigma_j^{(r-1)/2})(1 + \lambda^2 \mc O_{\lambda}) = 1 + \ell \cdot (r-1) \ZZ_{\ell}$.
We obtain the following decomposition:
\[
\mc T_{\ell, r} = T_{\ell, r}(\mu_{2 \ell}) \times (1 + \ell \cdot (r-1) \ZZ_{\ell}) \times \mc T'_{\ell, r}.
\]
Thus we are left with investigating $\mc U_{\ell}'$ and $\mc T'_{\ell, r}$. Recall that the $\ell$-adic logarithm
$\log : 1 + \lambda^2 \mc O_{\lambda} \to \lambda^2 \mc O_{\lambda}$ is an isomorphism of topological
groups, with the inverse given by the exponential function (see e.g. \cite[Example~IV.3.1.2, Theorem IV.6.4~(b)]{Silverman_AEC}). 
Hence we may identify $\mc U_{\ell}'$ with the set $\{ d \in \lambda^2 \mc O_{\lambda} : d + \ol d = 0 \}$
and $T_{\ell, r}'$ with the map $\lambda^2 \mc O_{\lambda} \to \lambda^2 \mc O_{\lambda}$ given by the formula $\sum_{j = 1}^{\ell - 1} n'(j) \cdot \sigma_j$.	
Note that $\lambda^2 \mc O_{\lambda}$ is a free $\ZZ_{\ell}$-module of rank $\ell - 1$. One easily proves that $\mc U_{\ell}' \subset \lambda^2 \mc O_{\lambda}$
is a free $\ZZ_{\ell}$-module of rank $\frac{\ell - 1}{2}$ with the basis given by
\[
\lambda^i - \ol{\lambda}^i \quad \textrm{ for } i = 2, \ldots, \frac{l+1}{2}.
\]
\begin{proof}[Proof of Theorem~\ref{thm:DC}]
	By Proposition~\ref{prop:Tlr_in_DJ} and the above discussion, $[\mc U_{\ell, r} : \mu_{2 \ell} \mc D_J] \le [\mc U_{\ell, r} : \mu_{2 \ell} \mc T_{\ell, r}] = [\mc U_{\ell}' : \mc T'_{\ell, r}]$.
	Consider the map
	\[
	T_{\ell, r}'' := T_{\ell, r}'|_{\mc U_{\ell}'} : \mc U_{\ell}' \to \mc U_{\ell}'.
	\]
	By the Smith normal theorem, $|\coker T_{\ell, r}''| = \ell^t$, where $t := \ord_{\ell}(\det T_{\ell, r}'')$.
	Let $C \subset (\ZZ/\ell)^{\times}$ be an arbitrary set of representatives for $(\ZZ/\ell)^{\times}/\{ 1, \ell-1 \}$.
	We will compute the matrix of $T_{\ell, r}'' \otimes \QQ_{\ell}$		
	in the basis $\{\zeta^i_{\ell} - \zeta^{-i}_{\ell} : i \in C \}$. For any $i \in C$ we have:
	\begin{align*}
		T_{\ell, r}''(\zeta^i_{\ell} - \zeta^{-i}_{\ell}) &= \sum_{j = 1}^{\ell - 1} n'(j) \cdot (\zeta_{\ell}^{ij} - \zeta^{-ij}_{\ell})
		= \sum_{j \in C} 2 \cdot n'(j) \cdot (\zeta_{\ell}^{ij} - \zeta_{\ell}^{-ij})\\
		&= \sum_{k \in C} 2 \cdot n'(i^{-1} k) \cdot (\zeta_{\ell}^k - \zeta_{\ell}^{-k}).
	\end{align*}
	(we used~ \eqref{eqn:n'_is_odd} and that $\ell-i \not \in C$ for $i \in C$).
	Hence, the matrix of $T_{\ell, r}''$ is given by $2 \cdot [n'(i^{-1} j)]_{i, j \in C}$.
	The proof follows from Lemma~\ref{lem:demjaenko_matrix} below.
\end{proof}
\begin{lemma} \label{lem:demjaenko_matrix}
	Let $\ell$ be an odd prime and $r$ a natural number non-divisible by $\ell$. Let $C \subset (\ZZ/\ell)^{\times}$ be an arbitrary
	set of representatives of $(\ZZ/\ell)^{\times}/\{ 1, \ell - 1 \}$. Define $n ' : (\ZZ/\ell)^{\times} \to \ZZ$ by the formula~\eqref{eqn:n_prim_eqn}
	and $c_{\ell, r}$ by~\eqref{eqn:c_lr_def}.
	Then:
	\[
	\det[n'(i \cdot j^{-1})]_{i, j \in C} = \frac{(-1)^{\frac{\ell - 1}{2}} \cdot h_{\ell}^-}{2 \ell} \cdot c_{\ell, r}.
	\]
\end{lemma}
\begin{proof}
	This is basically Hirabayashi's formula for the determinant of the Demjanenko matrix from~\cite{Hirabayashi_generalization_Maillet_Demyanenko}.
	We will use this formula in the form \cite[formula~(13)]{Kucera_Formulae_relative_class_nb}. In the notation of ibid. we take $(K, m, b) := (\QQ(\zeta_{\ell}), \ell, r)$.
	Note that then $K^+ = \QQ(\zeta_{\ell} + \zeta_{\ell}^{-1})$, $w_K = 2 \ell$, $Q_K = 2$. Hence:
	\[
	\det[n'(i \cdot j^{-1})]_{i, j \in C} = \frac{(-1)^{\frac{\ell - 1}{2}} \cdot h_{\ell}^-}{2 \ell} \cdot c_K(r)
	\]
	and we are left with proving that $c_K(r) = c_{\ell, r}$. Denote $R := r_{\ell}$ and $L := \frac{\ell - 1}{R}$.
	Then $r \equiv g^L \pmod{\ell}$ for a generator $g$ of $(\ZZ/\ell)^{\times}$. Moreover, the set
	of odd Dirichlet characters for $K$ is given by $X_- = \{ \chi_a : a = 1, \ldots, \ell - 1, \, 2 \nmid a \}$, where
	$\chi_a(g^i) := \zeta_{\ell - 1}^{a \cdot i}$. Therefore for any $\chi \in X_-$ we have $f_{\chi} = \ell$ and:
	\begin{align*}
		c_K(r) &= \prod_{\chi \in X_-} c_K^{\chi}(r) = \prod_{\chi \in X_-} (r - \chi(r))\\
		&= \prod_{\substack{a = 1\\2 \nmid a}}^{\ell - 1} (r - \zeta_{\ell - 1}^{L \cdot a})
		= \prod_{\substack{a = 1\\2 \nmid a}}^{\ell - 1} (r - \zeta_R^{a}).
	\end{align*}
	Thus, if $2 | R$, then:
	\begin{align*}
		c_K(r) &= \left(\prod_{\substack{a = 1\\2 \nmid a}}^R (r - \zeta_R^{a}) \right)^{L}
		= \left(
		\frac{\prod_{\substack{a = 1}}^R (r - \zeta_R^{a})}
		{\prod_{\substack{b = 1}}^{R/2} (r - \zeta_R^{2b})}
		\right)^{L}
		= \left(\frac{r^R - 1}{r^{R/2} - 1} \right)^{L}
		= (r^{R/2} + 1)^{L}.
	\end{align*}
	If $2 \nmid R$, then $2 | L$ and:
	\begin{align*}
		c_K(r) &= \left(\prod_{\substack{a = 1}}^R (r - \zeta_R^{a}) \right)^{L/2}
		= (r^R - 1)^{L/2}.
	\end{align*}
	This finishes the proof.
\end{proof}
Note that the Demjanenko matrix has been used previously in a related context
in order to study the Sato--Tate group of the Jacobians of Fermat curves,
see \cite{Fite_Gonzalez_Frobenius_distribution}.\\

We finish this section by computing the degree of the $\ell$-division field of the Jacobian
considered in Example~\ref{ex:intro}.
\begin{example} \label{ex:degree_division_field}
	Keep the notation of Example~\ref{ex:intro}. As explained in Example~\ref{ex:intro}, we have $\mc D_J = \mc U_{11, 8}$.
	Denote by $\mc U$ the reduction of $\mc U_{11, 8}$ modulo $11$. Then by Remark~\ref{rem:DJ_and_mu_2} and Lemma~\ref{lem:basis_of_su_n} we have:
	\begin{align*}
		[K(J[11]) : K] &=  |\GU_{7}^-(\mc O/\lambda^{10})^{\Gal(f)}_{\det \in \mc U}|\\
		&= |\Gal(f) \cap A_r| \cdot |\SU_{7}^-(\mc O/\lambda^{10})_1| \cdot |\mc U|\\
		&= \frac 12 8! \cdot \prod_{i = 1}^{9} |\SU_{7}^-(\mc O/\lambda^{i+1})_i|
		\cdot |\mc U|\\
		&= \frac 12 8! \cdot 11^{4 \cdot {7 \choose 2} + 5 \cdot ({8 \choose 2} - 1)} \cdot |\mc U|\\
		&= \frac 12 8! \cdot 11^{219} \cdot |\mc U|.
	\end{align*}
	Moreover, by the above discussion $\mc U_{11, 8} = \mu_{22} \times (1 + 11 \cdot \ZZ_{11}) \times \mc U_{11}'$,
	where $\mc U_{11}' \le 1 + \lambda^2 \mc O_{\lambda}$ is a $\ZZ_{11}$-lattice of rank $5$.
	This easily implies that $\mc U \cong \mu_{22} \times (\ZZ_{11}/11^8 \ZZ_{11})^5$ and $[K(J[11]) : K] = 8! \cdot 11^{260}$.
\end{example}
\section{The Hodge, Tate and Mumford-Tate conjectures} \label{sec:mt_conjecture}
In the last section we prove Corollary~\ref{cor:mt_conjecture}. We start by proving the Mumford--Tate conjecture 
for the Jacobians satisfying the assumptions of Theorem~\ref{thm:image_of_rho}.
Before the proof we recall some key notions. See e.g. \cite{Farfan_Survey_HTMTC} for a reference. The \bb{Mumford--Tate group} of $J$, denoted $\MT(J)$, is the largest algebraic subgroup of $\GSp(H_1(J, \QQ), \phi_J)$ fixing all the Hodge vectors. The \bb{Hodge group} of~$J$ is $\Hg(J) := (\MT(J) \cap \Sl(H_1(J, \QQ)))^{\circ}$. Let $\phi_J$ be the $\QQ$-linear symplectic pairing $H_1(J, \QQ) \times H_1(J, \QQ) \to \QQ$ defined
by the principal polarisation. Then $\Hg(J) \subset \Sp(H_1(J, \QQ), \phi_J)$ and
$\MT(J) = \GG_m \cdot \Hg(J)$. Write $V_{\ell} J := T_{\ell} J \otimes_{\ZZ_{\ell}} \QQ_{\ell}$ and define the \bb{$\ell$-adic Galois monodromy group} of $J$ as 
\[
H_{\ell} := \left( \ol{\rho_{J, \ell}(\Gal_K)}^{Zar} \cap \Sl(V_{\ell} J) \right)^{\circ}
\]
(from now on, $\Sl(V_{\ell} J)$ denotes the group of matrices $\bb A \in \Gl(V_{\ell} J)$ with $\det_{\ZZ_{\ell}} \bb A = 1$).
The Betti-\'{e}tale comparison isomorphism $H_1(J, \QQ) \otimes \QQ_{\ell} \cong V_{\ell} J$ leads to an inclusion $H_{\ell} \subset \Hg(J) \otimes_{\QQ} \QQ_{\ell}$. The Mumford--Tate conjecture is equivalent
to the equality of these groups. Moreover, the pairing $\phi_J$ can be upgraded to an $\mc O$-linear skew-Hermitian pairing $\Phi_J : H_1(J, \QQ) \times H_1(J, \QQ) \to \mc O$ 
similarly as $e_{\ell^{\infty}}$ in Lemma~\ref{lem:weil_general} (cf. \cite{Xue_Zarhin_Hodge_groups}). Under the comparison isomorphism
the pairing $\phi_J$ corresponds to the pairing $\mf e_{\ell^{\infty}}$ defined in Lemma~\ref{lem:weil_general}. 
Thus:
\[
\Hg(J) \otimes \QQ_{\ell} \subset \U(H_1(J, \QQ), \Phi_J) \otimes \QQ_{\ell}
\cong \U(V_{\ell} J, \mf e_{\ell^{\infty}}) = 
\U(V_{\ell} J, \bb e_{\ell^{\infty}})
\]
(in fact, the article \cite{Xue_Endomorphisms_superelliptic} shows that this inclusion is an equality if $\Gal(f) \cong S_r$ or~$A_r$).
Note that the unitary groups in question can be considered as algebraic groups over~$\QQ_{\ell}$ by
taking the Weil restriction from $\QQ_{\ell}(\zeta_{\ell} + \zeta_{\ell}^{-1})$ (cf. Remark~\ref{rmk:unitary_is_algebraic}) to $\QQ_{\ell}$.

On the other hand, note that
$\det_{\mc O_{\lambda}}(\U^{\varepsilon}_{r-1}(\mc O) \cap \Sl(T_{\ell} J)) = \mc U_{\ell}'$ (recall that $\mc U_{\ell}' := \{ d \in 1 + \lambda^2 \mc O_{\lambda} : d \cdot \ol d = 1 \}$). Hence by Theorem~\ref{thm:image_of_rho} and Remark~\ref{rem:DJ_and_mu_2}:
\[
\rho_{J, \ell}(\Gal_K) \cap \Sl(T_{\ell} J) = U(T_{\ell} J, \bb e_{\lambda^{\infty}})^{\Gal(f) \cap A_r}_{\det \in \mc D_J \cap \mc U_{\ell}'}.
\]
By Theorem~\ref{thm:DC}, the index of $\mc D_J \cap \mc U_{\ell}'$ in $\mc U_{\ell}'$ is finite.
Thus $\rho_{J, \ell}(\Gal_K) \cap \Sl(T_{\ell} J)$ has a finite index in $\U(T_{\ell} J, \bb e_{\lambda^{\infty}})$.
This easily implies that in $\Gl(V_{\ell} J)$:
\[
\ol{\rho_{J, \ell}(\Gal_K) \cap \Sl(T_{\ell} J)}^{Zar, \circ} = \ol{\U(T_{\ell} J, \bb e_{\lambda^{\infty}})}^{Zar, \circ}
= \U(V_{\ell} J, \bb e_{\lambda^{\infty}}),
\]
which ends the proof of the Mumford--Tate conjecture for $J$. As a consequence:
\begin{align*}
	\End(J) \otimes \QQ \cong \End_{\Hg(J)}(H^1(J, \QQ)) = \End_{U(H_1(J, \QQ), \Phi_J)}(H_1(J, \QQ)) = K.
\end{align*}
Thus, since $\mc O \subset \End(J)$, we conclude that $\End(J) = \mc O$.
Moreover, the Mumford--Tate conjecture holds also for all powers of $J$ by \cite{Commelin_MT_for_products}.
Finally, using \cite[Theorem 0]{Ribet_Hodge_classes_certain}, we conclude that the Hodge conjecture holds for~$J$
and all of its powers. The Tate conjecture follows from Hodge and Mumford--Tate conjectures, see e.g. \cite[Section~6]{Farfan_Survey_HTMTC}.
	
\bibliographystyle{apalike}	
\bibliography{bibliografia}

\newpage
\appendix
\section*{Appendix: Sage code} \label{appendix:sage}
The following Sage Math code verifies the identity from Lemma~\ref{lem:computation_of_commutator} for any given value of $n$.

\begin{python}
	n = 4 #change n to verify this for a different value
	N = floor((n-1)/2)
	#We define a noncommutative algebra
	#with generators AN, ... An-1, BN, ..., Bn-1
	generators = ''
	for i in range(N, n):
  generators += ('A'+str(i)+', ')
	for i in range(N, n):
  generators += ('B'+str(i))
  if i!=n-1:
      generators += ', '
	R = FreeAlgebra(QQ, 2*(n-N), names=generators)
	A = N*[0] + list(R.gens()[:n-N])
	B = N*[0] + list(R.gens()[n-N:])
	#We define a polynomial ring over R. 
	#The variable t plays the role of lambda.
	Rt.<t> = PolynomialRing(R)
	#Now we define the "matrices" with the given lambda-adic expansion.
	AA = 1 + sum(A[i]*t^i for i in range(n))
	BB = 1 + sum(B[i]*t^i for i in range(n))
	#We define now the commutator separately in the three cases
	if n
		commutator = 1 + t^(n-1)*(A[N]*B[N] - B[N]*A[N])
	if n
		commutator = 1 + t^(n-2)*(A[N]*B[N] - B[N]*A[N])
		+ t^(n-1)*(A[N]*B[N+1] - B[N+1]*A[N]
		+ A[N+1]*B[N] - B[N]*A[N+1])
	if n == 4:
		commutator = 1 + t^2*(A[1]*B[1] - B[1]*A[1])
		+ t^3*(A[1]*B[2] - B[2]*A[1] + A[2]*B[1] - B[1]*A[2]
		+ (B[1]*A[1] - A[1]*B[1])*(A[1] + B[1]))
	#We check now that AA*BB = commutator*BB*AA + O(t^n)
	print(list(AA*BB - commutator*BB*AA)[:n] == n*[0])
\end{python}
\end{document}